\documentclass[12pt,draftcls,onecolumn]{IEEEtran}
\usepackage{fancyhdr}

\def\baselinestretch{1.2}

\usepackage{dsfont}
\usepackage{amsfonts}
\usepackage{mathrsfs}
\usepackage{}
\usepackage{colortbl}
\usepackage{amssymb,amsmath,graphicx}
\usepackage{caption}

\usepackage{subfigure}
\usepackage{float}
\usepackage{multirow}

\newtheorem{Remark}{\bf Remark}[section]

\newtheorem{Problem}{\bf Problem}[section]

\newenvironment{Proof}{\noindent{\em Proof:\/}}{\hfill $\Box$\par}

\newtheorem{Theorem}{\bf Theorem}[section]

\newtheorem{Lemma}{\bf Lemma}[section]

\newtheorem{Proposition}{\bf Proposition}[section]

\newtheorem{Assumption}{\bf Assumption}[section]

\begin{document}
%
\title{Global Robust Practical Output Regulation for Nonlinear Systems in Output Feedback Form by Output-Based Event-Triggered Control}
%
%
%

\author{Wei~Liu,~\IEEEmembership{Member,~IEEE}~and~Jie~Huang,~\IEEEmembership{Fellow,~IEEE}
\thanks{This work was submitted to a journal on June 7, 2017, and was supported in part by the Research Grants Council of the Hong Kong Special Administration Region under grant No. 14200515, and in part by National Natural Science Foundation of China under Project 61633007.}
\thanks{Wei Liu and Jie Huang  are with the Department of Mechanical and Automation
Engineering, The Chinese University of Hong Kong, Shatin, N.T., Hong
Kong. E-mail: wliu@mae.cuhk.edu.hk, jhuang@mae.cuhk.edu.hk}
\thanks{Corresponding author: Jie Huang.} 
}

\maketitle

\thispagestyle{fancy}
\fancyhead{}
\lhead{}
\lfoot{}
\cfoot{}
\rfoot{}
\renewcommand{\headrulewidth}{0pt}
\renewcommand{\footrulewidth}{0pt}


\begin{abstract}
In this paper, we study the event-triggered global robust practical output regulation problem for a class of nonlinear systems in output feedback form with any relative degree. Our approach consists of the following three steps.  First, we design an internal model and an observer to form the so-called extended augmented system. Second, we  convert the original problem  into the event-triggered  global robust practical stabilization problem of the extended augmented system. Third, we  design an output-based event-triggered control law and a Zeno-free output-based event-triggered mechanism to solve the stabilization problem, which in turn leads to the solvability of the original problem.
Finally,  we apply our main result to the tracking problem of the controlled hyper-chaotic Lorenz systems.
\end{abstract}

%
\IEEEpeerreviewmaketitle



\section{Introduction}
\IEEEPARstart{T}{he} robust output regulation problem has been one of the fundamental and important control problems and has attracted extensive attention from the control community over the past few decades.
It aims to  design a feedback control law for an uncertain plant such that the output of the closed-loop system can asymptotically track some class of reference inputs while rejecting some class of  external disturbances. Here the reference inputs and external disturbances are both  generated by an autonomous differential equation called exosystem.
 The problem was first thoroughly studied for linear uncertain systems in, say,  \cite{Davision1976,Francis1977,Francis1976} during 1970s. The local robust output regulation problem for nonlinear systems was later studied in  \cite{Byrnes1997,Huang1995,Khalil1994} during 1990s. In \cite{Marconi2008}, the semi-global practical output regulation problem for nonlinear systems was studied. In \cite{Byrnes2004}, by using some nonlinear internal models, the semi-global robust output regulation problem for nonlinear systems was studied.
 In \cite{Huang2}, a general framework for tackling the global robust output regulation problem was established. The framework  consists of two steps: First, convert the global robust output regulation problem for a given plant to a global robust stabilization problem
for a well defined augmented system composed of the given plant and a specific dynamic compensator called internal model; Second, solve the global robust stabilization problem for the augmented system.  This framework has been successfully applied to several classes of typical
 nonlinear systems such as 
 lower triangular nonlinear systems \cite{ChenZHuang2004}, and output feedback nonlinear systems \cite{XuHuang2010c}.


Conventionally, a continuous-time control law is implemented in a digital platform through sampling the measured analog signal  with a fixed time period
\cite{Astrom1,Franklin1}. The sampling takes place without considering whether or not it is necessary, and, as a result, the samplings and control actions may be redundant, which leads to the  waste of the system resources.
In contrast, the event-triggered control is another way to implement  a continuous-time control law in a digital platform.
As reviewed in \cite{Heemels1}, the event-triggered control generates the samplings and the control actions only when the system state or output deviates away from a prescribed set or the performance index violates a specified level. Thus the event-triggered control strategy is able to reduce the control execution times and saving energy resources while maintaining the desired control performance.
A central technical issue with the event-triggered control is to exclude the Zeno behavior, that is, the execution times become arbitrarily close and converge to a finite  accumulation point \cite{Tabuada1}.

The event-triggered control method has now become an active research topic.
Various event-triggered control problems have been  widely studied for both linear and nonlinear systems. 
 For example,  reference \cite{Heemels1} first reviewed the event-triggered control approach and then designed a state-feedback event-triggered control law to solve the stabilization problem for a class of linear systems.
Reference \cite{Donkers2012} further designed an output-based event-triggered control law to guarantee  the $\mathcal{L}_{\infty}$-performance and the  closed-loop  stability for a class of linear systems.
Reference \cite{Tabuada1} solved the  stabilization problem for nonlinear systems by a state-based event-triggered control law under the ISS assumption with respect to the measurement state error.
In \cite{LiuT1,LiuT2},  the small gain theorem was used to solve the robust stabilization problem for nonlinear systems by a state-based event-triggered control law and an output-based event-triggered control law, respectively.
Reference \cite{Tallapragada2013} designed  a state-based event-triggered control law to solve the asymptotic tracking problem for  nonlinear systems  and  the tracking error was guaranteed to be  uniformly ultimately bounded.
In \cite{Abdelrahim2017},  the robust stabilization for nonlinear systems was studied by an output-based event triggered control law.
In particular, the global robust practical output regulation problem for nonlinear systems in normal form with unity relative degree was studied in our recent work \cite{LiuHuang2017b}.
 Some other  contributions relevant to this paper can be found in \cite{Dolk2016,XingL1}. 


In this paper,  we will further study  the global robust practical output regulation problem for  nonlinear systems in output feedback form with any relative degree  by an output-based event-triggered control law. The system to be studied is a generalization of the system studied in \cite{LiuHuang2017b}, which can be viewed as an output feedback system with unity relative degree. As will be seen in Remark \ref{RemarkProb1} and in the conclusion part of this paper. The problem to be studied in this paper is more complex than the one in \cite{LiuHuang2017b}. We need to establish some technical lemmas to
overcome some specific technical difficulties and develop a recursive approach to construct both the event-triggered  control law and
the event-triggered mechanism.


Compared with some other existing event-triggered control problems for nonlinear systems such as stabilization problem \cite{LiuT1,Tabuada1} 
 and tracking problem \cite{Tallapragada2013,XingL1},
our problem poses at least three specific challenges.  First, the control objective of our problem is to achieve not only the asymptotical tracking  and but also the disturbance rejection for a class of uncertain nonlinear systems in output feedback form with any relative degree. Second, the uncertain parameter vector of the plant is allowed to belong to an arbitrarily large prescribed compact set.
Third, our control law is a dynamic output feedback control law which contains not only an internal model but also an observer.  Thus, we need
to sample not only the output of the plant but also the states of the internal model and the observer. As a result, the stability analysis of the closed-loop system is more complicated than the static state or static output feedback case. We have managed to overcome these challenges by integrating the internal model approach and the observer-based approach, and shown that our design can exclude the Zeno phenomenon.

The rest of this paper is organized as follows. In Section \ref{PF}, we give the problem formulation
and some preliminaries. In Section \ref{MR},  we solve the problem  by a recursively designed event-triggered output feedback control law together with an output-based event-triggered mechanism. In Section \ref{Example}, a simulation example is given to to illustrate the design. In Section \ref{Conclusion}, we give some concluding remarks. Finally,  one technical lemma and some proofs are given in appendix section.

{\bf Notation.} For any column vectors $a_i$, $i=1,...,s$, denote $\mbox{col}(a_1,...,a_s)=[a_1^T,...,a_s^T]^T$.
The set of all nonnegative integers is denoted by $\mathbb{N}$.
 The base of the natural logarithm is denoted  by $\textbf{e}$.
The maximum eigenvalue and the minimum eigenvalue of a symmetric real matrix $A$ are denoted by  $\lambda_{\max}(A)$ and $\lambda_{\min}(A)$, respectively.
In this paper, for simplicity, we use $x$ to denote $x(t)$ when no ambiguity occurs.

\section{Problem formulation and preliminaries}\label{PF}
Consider a class of nonlinear  systems in output feedback form  as follows
\begin{equation}\label{system1}
\begin{split}
  \dot{z}  =& ~f (z ,y ,v,w)\\
  \dot{x}_{i}=&~g_{i}(z,y,v,w)+x_{i+1},~i=1,\cdots,r-1\\
  \dot{x}_{r}  =&~ g_{r} (z ,y,v,w) + b(w)u \\
  y =&~x_{1}\\
\end{split}
\end{equation}
where $z\in \mathbb{R}^{n_{z}}$ and $x=\mbox{col}(x_{1},\cdots,x_{r})\in \mathbb{R}^{r}$ are the states,  $y \in \mathbb{R}$ is the output, $u \in \mathbb{R}$ is the input,
$w\in\mathbb{R}^{n_{w}}$ is an uncertain constant vector, and $v(t)\in\mathbb{R}^{n_{v}}$ is an exogenous signal representing both reference input  and external disturbance.
It is assumed that $v(t)$ is  generated by  a linear exosystem as follows
\begin{equation}\label{exosystem1}
\begin{split}
  \dot{v} = Sv,\ \ y_{0}=q(v,w).
\end{split}
\end{equation}
Define the regulated error output as $e=y-y_{0}$. We assume that all functions in (\ref{system1}) and (\ref{exosystem1}) are sufficiently smooth, and satisfy $q(0,w)=0$, $f(0,0,0,w)=0$, $g_{i} (0,0,0,w)=0$ with $i=1,\cdots,r$, and $b(w)>0$ for all $w\in \mathbb{R}^{n_{w}}$.

Consider a control law of the following form
\begin{equation}\label{u1}
\begin{split}
u (t)&=\hat{f} (\eta(t_{k}),\hat{\xi}(t_{k}),e(t_{k}))\\
\dot{\eta} (t)&=\hat{g}(\eta(t),\eta(t_{k}),\hat{\xi}(t_{k}),e(t_{k}))\\
\dot{\hat{\xi}} (t)&=\hat{l}(\hat{\xi}(t),\hat{\xi}(t_{k}),e(t_{k})),~t\in[t_{k},t_{k+1}),~k\in\mathbb{S}\\
\end{split}
\end{equation}
where   $\hat{f} (\cdot)$, $\hat{g}(\cdot)$ and $\hat{l}(\cdot)$ are some nonlinear functions,  $\eta \in \mathbb{R}^{s}$ and $\hat{\xi}\in \mathbb{R}^{r}$ are the states of the internal model and observer which will be described in detail later in this section, $\mathbb{S}\subseteq\mathbb{N}$ denotes the set of triggering times, and $\{t_{k}\}_{k\in\mathbb{S}}$ denotes the triggering time sequence with $t_{0} =0$, which  is generated by an event-triggered mechanism of the following form
\begin{equation}\label{trigger1}
\begin{split}
t_{k+1} =\inf\{t>t_{k}~|~\hat{h} (\tilde{\eta}(t),\tilde{\xi}(t),\tilde{e}(t),\hat{\xi}(t),e(t))\geq\delta\}
\end{split}
\end{equation}
where    $\hat{h} (\cdot)$ is some nonlinear function, $\delta>0$ is some constant, and
 \begin{equation}\label{tildee1}
\begin{split}
 &\tilde{e}(t)=e(t_{k} )-e(t)\\
 &\tilde{\eta}(t)=\eta(t_{k})-\eta(t)\\
  &\tilde{\xi}(t)=\hat{\xi}(t_{k})-\hat{\xi}(t),~\forall t\in[t_{k},t_{k+1}),~k\in\mathbb{S}.\\
\end{split}
\end{equation}
\begin{figure}[H]
\centering
\includegraphics[scale=0.8]{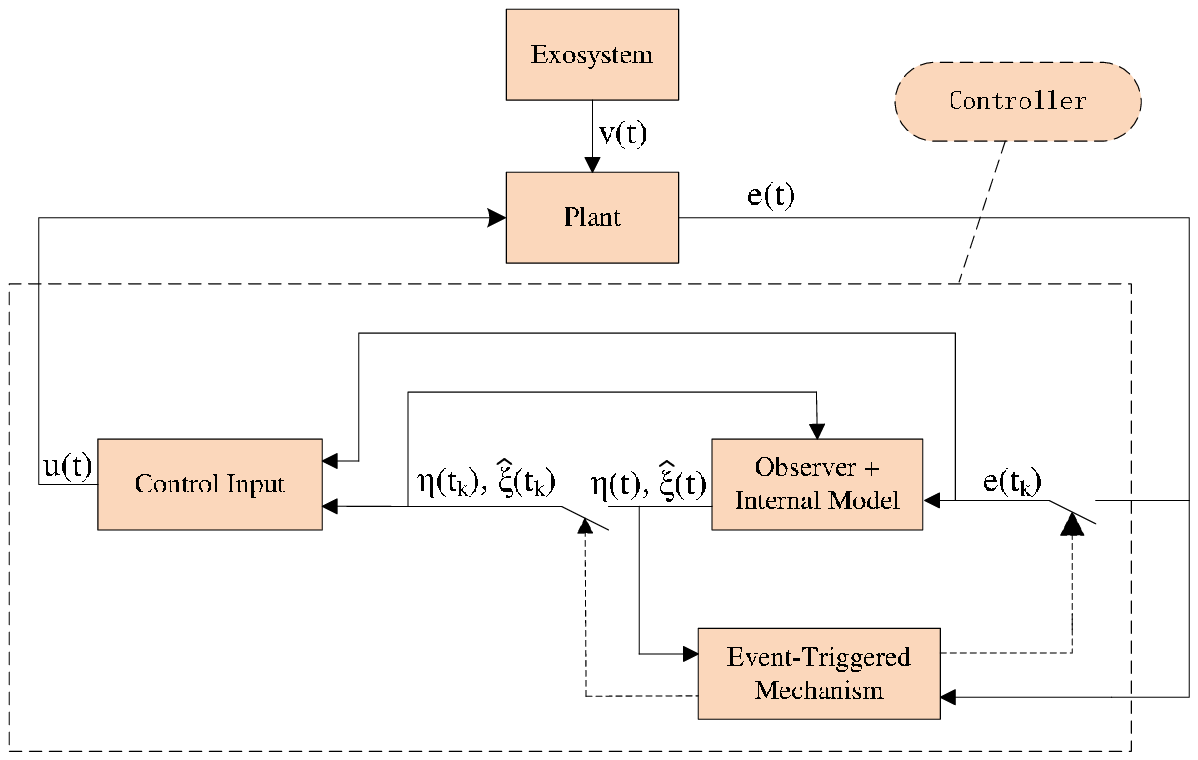}
\caption{Event-triggered  control schematic.} \label{diagram1}
\end{figure}
The closed-loop system under our event-triggered  control scheme is given in Figure \ref{diagram1}.

Now we describe our problems as follows.
%
\begin{Problem}\label{Problem1}
 Given the plant (\ref{system1}), the exosystem (\ref{exosystem1}),    some compact subsets
$\mathbb{V}\subset\mathbb{R}^{n_{v}}$ and $\mathbb{W}\subset\mathbb{R}^{n_{w}}$ containing the origin, 
 and any $\epsilon>0$, design a control law of the form (\ref{u1}) and an event-triggered mechanism of the form \eqref{trigger1}   such that,
  for any $v\in\mathbb{V}$, $w\in\mathbb{W}$, and any initial states $z (0)$, $x (0)$, $\eta (0)$, $\hat{\xi}(0)$,
\begin{enumerate}
  \item the trajectory of the closed-loop system  exists and is bounded for all $t\geq0$;
  \item $\lim_{t \to \infty}\sup|e(t)|\leq \epsilon$.
\end{enumerate}
\end{Problem}

\begin{Remark}\label{RemarkProb1}
As in \cite{LiuHuang2017b}, Problem \ref{Problem1} is called as the event-triggered  global robust practical output regulation problem, and a control law that solves Problem \ref{Problem1} is called as a practical solution to the global robust output regulation problem.
The problem in \cite{LiuHuang2017b} can be viewed as a special case of the problem here since the system in \cite{LiuHuang2017b} is a special case of the system \eqref{system1} here with $r=1$.  Like in \cite{LiuHuang2017b}, the first step of our design is to employ the internal model approach to convert the problem into a stabilization problem of a well defined augmented system. However, unlike in
\cite{LiuHuang2017b} where the augmented system  can be stabilized by a static output feedback control law, here we need to further employ a hybrid practical observer to estimate the partial state of the augmented system, which leads to a hybrid extended augmented system.
Then an event-triggered static control law will be further designed to stabilize the extended augmented system.  Moreover, since the system  \cite{LiuHuang2017b} has unity relative degree, the design of the control law and the event-triggered mechanism does not involve a recursive procedure.
In contrast, here both of our event-triggered control law of the form \eqref{u1} and our event-triggered mechanism of the form \eqref{trigger1} have to be designed recursively,
which makes the stability analysis of the closed-loop system much more complicated.

\end{Remark}


To solve our problem, 
we first introduce some standard assumptions which can also be found in \cite{LiuHuang2017b,XuHuang2010c}.

\begin{Assumption}\label{Ass2.1}
The exosystem is neutrally stable, i.e., all the eigenvalues of $S$ are semi-simple with zero real parts.
\end{Assumption}

Under Assumption \ref{Ass2.1}, for any $v(0)\in\mathbb{V}_{0}$ with $\mathbb{V}_{0}$ being some known compact set, there always exists another compact set $\mathbb{V}$ such that $v(t)\in\mathbb{V}$ for all $t\geq0$.

%

\begin{Assumption}\label{Ass2.3}
There exists a globally defined smooth function $\textbf{z} :\mathbb{R}^{n_{v}}\times\mathbb{R}^{n_{w}}\mapsto\mathbb{R}^{n}$ with $\textbf{z} (0,w)=0$ such that
\begin{equation}\label{Z}
\begin{split}
  \dfrac{\partial\textbf{z} (v,w)}{\partial v}Sv=f (\textbf{z} (v,w),q(v,w),v,w)
\end{split}
\end{equation}
for all $(v,w)\in\mathbb{R}^{n_{v}}\times\mathbb{R}^{n_{w}}$.
\end{Assumption}

Under Assumption \ref{Ass2.3}, let
\begin{equation*}
  \begin{aligned}
    &\textbf{x}_{1}(v,w)=q(v,w)\\
    &\textbf{x}_{i}(v,w)=\frac{\partial \textbf{x}_{i-1}(v,w)}{\partial v}Sv-g_{i-1}(\textbf{z}(v,w),q(v,w),v,w),~i=2,\cdots,r\\
    &\textbf{x}(v,w)=\mbox{col}(\textbf{x}_{1}(v,w),\cdots,\textbf{x}_{r}(v,w))\\
    &\textbf{u}(v,w)=b^{-1}(w)\big(\frac{\partial \textbf{x}_{r}(v,w)}{\partial v}Sv-g_{r}(\textbf{z}(v,w),q(v,w),v,w)\big).\\
  \end{aligned}
\end{equation*}
Then the solution to the regulator equations associated with  (\ref{system1})
and (\ref{exosystem1}) is given by $\mbox{col}(\textbf{z}(v,w)$, $\textbf{x}(v,w),\textbf{u}(v,w))$ \cite{Byrnes1997, Huang1}.


\begin{Assumption}\label{Ass2.4}
 The function $\textbf{u} (v,w)$ is a polynomial in $v$ with coefficients depending on $w$.
\end{Assumption}

Under Assumption \ref{Ass2.4}, there exists an integer $s$ such that, for all trajectories $v(t)$ of the exosystem \eqref{exosystem1} and all $w\in\mathbb{R}^{n_{w}}$, $\textbf{u} (v,w)$ satisfies
$\dfrac{d^{s }\textbf{u}(v(t),w) }{dt^{s }}=\varrho_{1}\textbf{u}(v(t),w) +\varrho_{2}\dfrac{d\textbf{u}(v(t),w) }{dt}+\cdots+\varrho_{s }\dfrac{d^{s -1}\textbf{u}(v(t),w) }{dt^{s -1}}$
where  $\varrho_{1},\cdots,\varrho_{s}$ are some real scalars such that all the roots of the polynomial $P (\lambda)=\lambda^{s }-\varrho_{1}-\varrho_{2}\lambda-\cdots-\varrho_{s }\lambda^{s -1}$  are distinct with zero real part \cite{Huang1}.
Let 
\begin{equation*}
\begin{split}
 \Phi =\left[
            \begin{array}{cccc}
              0 & 1 & \cdots & 0 \\
              \vdots & \vdots & \ddots & \vdots \\
              0 & 0 & \cdots & 1 \\
              \varrho_{1} & \varrho_{2} & \cdots & \varrho_{s } \\
            \end{array}
          \right],~ \Gamma =\left[
                                 \begin{array}{c}
                                   1 \\
                                   0 \\
                                   \vdots \\
                                   0 \\
                                 \end{array}
                               \right]^{T}.
\end{split}
\end{equation*}
Then the Sylvester equation $T \Phi -M T =N \Gamma$
has a unique nonsingular solution $T$ for any given  controllable pair $(M,N)$ with $M \in \mathbb{R}^{s \times s} $ a Hurwitz matrix and $N  \in \mathbb{R}^{s \times 1} $ a column vector \cite{Huang1}.
We further let $\Psi =\Gamma T ^{-1}$ and $\theta (v,w)=T \mbox{col}(\textbf{u} (v,w),\dot{\textbf{u}} (v,w),\cdots,\textbf{u} ^{(s -1)}(v,w))$.
Then we have 
$\dot{\theta}(v,w)=(M+N\Psi)\theta(v,w)$ and $\textbf{u} (v, w) = \Psi \theta (v, w)$.
Moreover,  define the following dynamic compensator
\begin{equation}\label{doteta1}
\begin{split}
  \dot{\eta} =M \eta +N u
\end{split}
\end{equation}
which is called as a linear internal model of \eqref{system1} like in \cite{Byrnes1997, Huang1, nik98}.

Motivated by \cite{XuHuang2010c}, we perform the following coordinate and input transformation on the system \eqref{system1} and the internal model \eqref{doteta1}
\begin{equation}\label{transformation1}
\begin{split}
  \bar{z} &=z -\textbf{z} (v,w),~\bar{x} =x -\textbf{x} (v,w),\\
  \bar{\eta}& =\eta -\theta (v,w)-C\bar{x}\\
\bar{u} &=u -\Psi \eta(t_{k}),~ t\in[t_{k},t_{k+1}),k\in\mathbb{S}\\
\end{split}
\end{equation}
where $C=[c_{1}~c_{2}~\cdots~c_{r}]$ with $c_{r}=b^{-1}(w)N$, $c_{i-1}=Mc_{i}$ for $i=2,\cdots,r$.
Then we obtain an augmented system as follows
\begin{equation}\label{system2}
\begin{split}
  \dot{\bar{z}} =&\bar{f} (\bar{z} ,e,v,w)\\
  \dot{\bar{\eta}} =&M \bar{\eta} +M c_{1}e-\sum_{i=1}^{r}c_{i}\bar{g}_{i}(\bar{z},e,v,w)\\
  \dot{\bar{x}} =&A_{d}\bar{x}+b(w)B\Psi\bar{\eta}+\bar{g}(\bar{z},e,v,w)+b(w)B\bar{u}+b(w)B\Psi\tilde{\eta}\\
\end{split}
\end{equation}
where $\bar{f} (\bar{z} ,e ,v,w)=f (\bar{z} +\textbf{z}(v,w) ,e +q(v,w),v,w)-$ $f (\textbf{z}(v,w) ,q(v,w),v,w)$, $\bar{g}(\bar{z} ,e ,v,w)=\mbox{col}(\bar{g}_{1} (\bar{z} ,e ,v,w)$, $\cdots,\bar{g}_{r} (\bar{z} ,e ,v,w))$, $\bar{g}_{i} (\bar{z} ,e ,v,w)=g_{i} (\bar{z} +\textbf{z}(v,w) ,e +q(v,w),v,w)-g_{i} (\textbf{z}(v,w)$, $q(v,w),v,w)$ for $i=1,\cdots,r$,
  and
\begin{equation*}
\begin{split}
A_{d}\!=\!\!\left[
    \begin{array}{ccccc}
      0 & 1 & 0 & \cdots & 0 \\
      0 & 0 & 1 & \cdots & 0 \\
      \vdots & \vdots & \vdots & \ddots & \vdots \\
      0 & 0 & 0 & \cdots & 1 \\
      d_{r} & d_{r-1} & d_{r-2} & \cdots & d_{1} \\
    \end{array}
  \right]_{r\times r}\!\!\!\!\!\!\!\!\!,~~~B\!=\!\!\left[
               \begin{array}{c}
                 0 \\
                 0 \\
                 \vdots \\
                 0 \\
                 1 \\
               \end{array}
             \right]_{r\times1}
\end{split}
\end{equation*}
with $d_{i}=b(w)\Psi c_{r+1-i}$ for $i=1\cdots,r$. Clearly, $\bar{f} (0,0,v,w)$ $=0$ and $\bar{g} (0,0,v,w)=0$ for any $v\in\mathbb{R}^{n_{v}}$ and $ w\in\mathbb{R}^{n_{w}}$.
Then, as in \cite{XuHuang2010c}, we further perform another coordinate transformation on the $\bar{x}$-subsystem as follows
  \begin{equation}\label{transformation2}
\begin{split}
\xi=b^{-1}(w)U_{d}\bar{x}
\end{split}
\end{equation}
where
\begin{equation*}
\begin{split}
U_{d}=\left[
    \begin{array}{ccccc}
      1 & 0 &  \cdots & 0 & 0\\
      -d_{1} & 1 & \cdots & 0 & 0 \\
      \vdots & \vdots & \ddots & \vdots & \vdots \\
      -d_{r-2} & -d_{r-3} & \cdots & 1 & 0 \\
      -d_{r-1} & -d_{r-2} & \cdots &   -d_{1} &1 \\
    \end{array}
  \right]_{r\times r}\!\!\!\!\!\!.
\end{split}
\end{equation*}
Then we have
\begin{equation}\label{system3}
\begin{split}
  \dot{\bar{z}} =&\bar{f} (\bar{z} ,e,v,w)\\
  \dot{\bar{\eta}} =&M \bar{\eta} +M c_{1}e-\sum_{i=1}^{r}c_{i}\bar{g}_{i}(\bar{z},e,v,w)\\
  \dot{\xi} =&A_{c}\xi+B\Psi\bar{\eta}+G(\bar{z},e,v,w)+B\bar{u}+B\Psi\tilde{\eta}\\
\end{split}
\end{equation}
where
 $G(\bar{z},e,v,w)=\mbox{col}(G_{1}(\bar{z},e,v,w),\cdots,G_{r}(\bar{z},e,v,w))$,
$G_{1}(\bar{z},e,v,w)=b^{-1}(w)(d_{1}e+\bar{g}_{1}(\bar{z},e,v,w))$, $G_{i}(\bar{z},e,v,w)$ $=b^{-1}(w)(d_{i}e-\sum_{j=1}^{i-1}d_{i-j}\bar{g}_{j}(\bar{z},e,v,w)+\bar{g}_{i}(\bar{z},e,v,w))$ for $i=2,\cdots,r$, and $A_{c}=\left[
                             \begin{array}{cc}
                               \textbf{0}_{(r-1)\times1} & I_{r-1} \\
                               0 & \textbf{0}_{1\times(r-1)} \\
                             \end{array}
                           \right]
$.
\begin{Remark}\label{RemarkTransfornation1}
The transformation \eqref{transformation1} is to convert the robust output regulation problem of \eqref{system1} into the robust stabilization problem of the augmented system \eqref{system2}. This step follows from the general framework for handling the robust output regulation problem described in \cite{Huang2}.
The second step is to globally stabilize \eqref{system2}. But since \eqref{system2} does not take any standard form that is amenable to some known stabilization technique, we further perform
the coordinate transformation \eqref{transformation2} to convert \eqref{system2} to a more standard lower triangular form \eqref{system3} with $\bar{z},\bar{\eta}$ as dynamic uncertainty and $B(\bar{u}+\Psi\tilde{\eta})$ as input. Here the transformation \eqref{transformation2} follows from the same one in \cite{XuHuang2010c}.
\end{Remark}

   The stabilization problem of systems of the form \eqref{system3} without the term $B\Psi\tilde{\eta}$
   has been well studied in the literature \cite{ChenZHuang2015, XuHuang2010c} by continuous-time control laws.
  Here, we will further consider
 stabilizing the system \eqref{system3} by an event-triggered control law.
 For this purpose,
mimicking the approach in \cite{ChenZHuang2015, XuHuang2010c}, we attach a dynamic compensator for \eqref{system3} as follows:
\begin{equation}\label{observer1}
\begin{split}
\dot{\hat{\xi}}=A_{o}\hat{\xi}+\lambda e(t_{k})+B\bar{u}
\end{split}
\end{equation}
where $\lambda=\mbox{col}(\lambda_{1},\cdots,\lambda_{r})$ and $\lambda_{1},\cdots,\lambda_{r}$ are chosen such that the matrix 
  \begin{equation*}
\begin{split}
A_{o}=\left[
    \begin{array}{ccccc}
      -\lambda_{1} & 1 & 0 & \cdots & 0 \\
       -\lambda_{2} & 0 & 1 & \cdots & 0 \\
      \vdots & \vdots & \vdots & \ddots & \vdots \\
       -\lambda_{r-1} & 0 & 0 & \cdots & 1 \\
       -\lambda_{r} & 0 & 0 & \cdots & 0 \\
    \end{array}
  \right]_{r\times r}\\
\end{split}
\end{equation*}
is Hurwitz. 

\begin{Remark}\label{RemarkObserver}
Like in \cite{ChenZHuang2015, XuHuang2010c}, the dynamic compensator (\ref{observer1}) is to estimate the state $\xi$ of the system  \eqref{system3}, which cannot be used for control  since it  relies on $v$ and $w$ according to \eqref{transformation1} and \eqref{transformation2}. However, unlike in
\cite{ChenZHuang2015, XuHuang2010c}, the dynamic compensator \eqref{observer1} is not continuous as it also relies on
the discrete quantity $e(t_{k})$. Nevertheless, it will be shown in the next section that the quantity $\|\xi-\hat{\xi}\|$ can still be made to approach an arbitrarily small value as $t$ tends to infinity. Thus, the system  (\ref{observer1}) can be viewed as a practical observer of the system  \eqref{system3}.
\end{Remark}

Now define the observation error as $\bar{\xi}=\xi-\hat{\xi}$. Then, from \eqref{tildee1}, \eqref{system3} and \eqref{observer1}, we have
\begin{equation}\label{observationerror1}
\begin{split}
\dot{\bar{\xi}}=&~A_{o}\bar{\xi}+\lambda(b^{-1}(w)-1)e+B\Psi\bar{\eta}+G(\bar{z},e,v,w)+B\Psi\tilde{\eta}-\lambda\tilde{e}.
\end{split}
\end{equation}
Attach \eqref{observationerror1} to \eqref{system3} and replace the state variable $\xi$ by $\mbox{col}(e,\hat{\xi}_{2},\cdots,\hat{\xi}_{r})$. Then we obtain the following system:
\begin{equation}\label{system4}
\begin{split}
&\dot{\bar{z}} =\bar{f} (\bar{z} ,e,v,w)\\
&\dot{\bar{\eta}} =M \bar{\eta} +M c_{1}e-\sum_{i=1}^{r}c_{i}\bar{g}_{i}(\bar{z},e,v,w)\\
&\dot{\bar{\xi}}=A_{o}\bar{\xi}+\lambda(b^{-1}(w)-1)e+B\Psi\bar{\eta}+G(\bar{z},e,v,w)+B\Psi\tilde{\eta}-\lambda\tilde{e}\\
&\dot{e}=b(w)\bar{\xi}_{2}+b(w)\hat{\xi}_{2}+b(w)G_{1}(\bar{z},e,v,w)\\
&\dot{\hat{\xi}}_{i}=\hat{\xi}_{i+1}+\lambda_{i}(e-\hat{\xi}_{1})+\lambda_{i}\tilde{e},~i=2,\cdots,r-1\\
&\dot{\hat{\xi}}_{r}=\bar{u}+\lambda_{r}(e-\hat{\xi}_{1})+\lambda_{r}\tilde{e}.\\
\end{split}
\end{equation}
\begin{Remark}\label{RemarkSystem4}
We call  (\ref{system4}) the extended augmented system of the given plant \eqref{system1}. It is noted that the derivation of the extended augmented system
(\ref{system4}) is different from that of the extended augmented system in \cite{XuHuang2010c} in two ways.
First the transformation  \eqref{transformation1} is different from the corresponding transformation in \cite{XuHuang2010c} in that
we have replaced the continuous function  $\eta(t)$ by the piecewise constant function $\eta(t_{k}),~k\in\mathbb{S}$. Second, the observer \eqref{observer1} is obtained from the corresponding observer in \cite{XuHuang2010c} by
replacing the continuous function $e(t)$ with the piecewise constant function $e(t_{k}),~k\in\mathbb{S}$. As a result, the extended augmented system
(\ref{system4}) here is also quite different from  the extended augmented system in \cite{XuHuang2010c} in that
it contains some additional terms $\tilde{\eta}$ and $\tilde{e}$, and thus is a hybrid system.
\end{Remark}

To stabilize (\ref{system4}) by an event triggered control law, we consider a piecewise constant control law as follows:
  \begin{equation}\label{baru1}
\begin{split}
\bar{u} (t)&=\check{f} (\hat{\xi}(t_{k}),e(t_{k})),~ t\in[t_{k},t_{k+1}),~k\in\mathbb{S}\\
\end{split}
\end{equation}
 where $\check{f} (\cdot)$ is a globally defined sufficiently smooth function vanishing at the origin.
Denote the state of the closed-loop system composed of (\ref{system4}) and (\ref{baru1}) under the triggering mechanism \eqref{trigger1} by
$\bar{x}_c = \mbox{col } (\bar{z},\bar{\eta},\bar{\xi},e,\hat{\xi}_{2},\cdots$, $\hat{\xi}_{r})$.
Then we first establish the following proposition.
\begin{Proposition}\label{Proposition1}
Under Assumptions \ref{Ass2.1}-\ref{Ass2.4}, for any $\epsilon>0$, any known compact sets $\mathbb{V}\in\mathbb{R}^{n_{v}}$ and $\mathbb{W}\in\mathbb{R}^{w}$,  if we can find a control law of the form (\ref{baru1}) and an event-triggered mechanism of the form \eqref{trigger1} such that, for any $\bar{x}_c (0)$, any $v(t)\in\mathbb{V}$ and any $w\in\mathbb{W}$, the solution $\bar{x}_{c}(t)$ exists and is bounded for all $t\in[0,\infty)$, and satisfies
\begin{equation}\label{pcgs1}
\lim_{t \to \infty}\sup\|\bar{x}_c(t)\|\leq \epsilon,
\end{equation}
then Problem \ref{Problem1} 
 for the original system \eqref{system1} is solvable by the following control law 
 \begin{equation}\label{u2}
\begin{split}
u (t)&=\check{f} (\hat{\xi}(t_{k}),e(t_{k}))+\Psi \eta (t_{k})\\
\dot{\hat{\xi}}(t)&=A_{o}\hat{\xi}(t)+\lambda e(t_{k})+B\check{f} (\hat{\xi}(t_{k}),e(t_{k}))\\
\dot{\eta} (t)&=M \eta (t)+N u (t),~t\in[t_{k},t_{k+1}),~k\in\mathbb{S}\\
\end{split}
\end{equation}
under the event-triggered mechanism \eqref{trigger1}.
\end{Proposition}

To not distract the attention of the readers,  we put the proof of Proposition \ref{Proposition1} in Appendix \ref{appendix1}.
As in \cite{LiuHuang2017b}, we call the problem of designing a control law of the form (\ref{baru1}) and an event-triggered mechanism of the form \eqref{trigger1}  to achieve (\ref{pcgs1}) as the  event-triggered global robust practical stabilization problem for   \eqref{system4}.

\section{Main Result}\label{MR}

 By Proposition \ref{Proposition1}, to solve the global robust practical output regulation problem for the system \eqref{system1}, it suffices to solve the global robust practical stabilization for the extended augmented system \eqref{system4}. In this section,
we will first establish a technical lemma which is given in Appendix \ref{appendix2}, and then develop a  recursive procedure to construct a piecewise constant output feedback control law and an output-based event-triggered mechanism. Finally, we will show that the solution of the closed-loop system exists for all $t\in[0,\infty)$, thus excluding the Zeno behavior.

For this purpose,  we introduce one more assumption.
 \begin{Assumption}\label{Ass3.1}
 For any compact subset $\Omega \subset \mathbb{R}^{n_{v}}\times\mathbb{R}^{n_{w}}$, there exists a $\mathcal{C}^{1}$ function $V_{0}(\bar{z} )$ such that, for  any $\mbox{col}(v,w)\in \Omega$, and any
$\bar{z}$ and $e$,
   \begin{equation}\label{Vz0}
\begin{split}
\underline{\alpha}_{0}(\|\bar{z} \|)\leq V_{0}(\bar{z} )
\leq\bar{\alpha}_{0}(\|\bar{z} \|)
\end{split}
\end{equation}
 \begin{equation}\label{dotVz0}
\begin{split}
\frac{\partial V_{0}(\bar{z})}{\partial \bar{z}} \bar{f} (\bar{z} ,e,v,w) \leq
-\alpha_{0}(\|\bar{z} \|)+\gamma_{0}(e )
\end{split}
\end{equation}
where $\gamma_{0}(\cdot)$ is a known smooth  positive definite function,  $\underline{\alpha}_{0}(\cdot)$, $\bar{\alpha}_{0}(\cdot)$ and $\alpha_{0}(\cdot)$ are  some known class $\mathcal{K}_{\infty}$ functions with $\alpha_{0}(\cdot)$ satisfying $\lim_{s\rightarrow0^{+}}\sup(s^{2}/\alpha_{0}(s))<\infty$.
 \end{Assumption}
 \begin{Remark}\label{RemarkAss3.1}
  Assumption \ref{Ass3.1} is  a standard assumption for nonlinear stabilization problem and has also been used in \cite{LiuHuang2017b,XuHuang2010c}.  Under this Assumption, the subsystem
$\dot{\bar{z}} =\bar{f} (\bar{z} ,e,v,w)$ is  input-to-state stable (ISS) with $e $ as the input \cite{Sontag1}.
\end{Remark}


Motivated by \cite{XuHuang2010c}, our control law will be recursively constructed as follows:
\begin{equation}\label{transformation3}
\begin{split}
&\check{\xi}_{1}(t)=e(t)\\
&\check{\xi}_{i+1}(t)=\hat{\xi}_{i+1}(t)-\vartheta_{i}(\check{\xi}_{i}(t)),~i=1,\cdots,r-1\\
\end{split}
\end{equation}
and
 \begin{equation}\label{baru2}
\begin{split}
&\bar{u}(t)=\vartheta_{r}(\check{\xi}_{r}(t_{k})),~ t\in[t_{k},t_{k+1}),~k\in\mathbb{S}\\
\end{split}
\end{equation}
where $\vartheta_{1}(\cdot),\cdots,\vartheta_{r}(\cdot)$ are some smooth functions to be constructed in Lemma \ref{Lemma1} below.   Let $X_{0}=\mbox{col}(\bar{z},\bar{\eta},\bar{\xi})$, $X_{i}=\mbox{col}(X_{i-1},\check{\xi}_{i})$ for $i=1,\cdots,r$, $\chi=\mbox{col}(\eta,e)$, $\tilde{\chi}=\mbox{col}(\tilde{\eta},\tilde{e})$, $\mu=\mbox{col}(v,w)$ and $\tilde{\vartheta}_{r}(t)=\vartheta_{r}(\check{\xi}_{r}(t_{k}))-\vartheta_{r}(\check{\xi}_{r}(t))$ for $t\in[t_{k},t_{k+1})$ with $k\in\mathbb{S}$.

Under the  transformation \eqref{transformation3}, the closed-loop system composed of \eqref{system4} and \eqref{baru2} can be put into the following form:
\begin{equation}\label{system5}
\begin{split}
&\dot{\bar{z}} =\bar{f} (\bar{z} ,e,\mu)\\
&\dot{\bar{\eta}} =M \bar{\eta} +\bar{f}_{\bar{\eta}}(\bar{z},e,\mu)\\
&\dot{\bar{\xi}}=A_{o}\bar{\xi}+B\Psi\bar{\eta}+\bar{f}_{\bar{\xi}}(\bar{z},e,\tilde{\chi},\mu)\\
&\dot{e}=b(w)(\check{\xi}_{2}+\vartheta_{1})+h_{1}(X_{0},e,\mu)\\
&\dot{\check{\xi}}_{i}=\check{\xi}_{i+1}+\vartheta_{i}+h_{i}(X_{i-1},\check{\xi}_{i},\tilde{\chi},\mu),~i=2,\cdots,r-1\\
&\dot{\check{\xi}}_{r}=\vartheta_{r}+\tilde{\vartheta}_{r}+h_{r}(X_{r-1},\check{\xi}_{r},\tilde{\chi},\mu)\\
\end{split}
\end{equation}
where
$\bar{f}_{\bar{\eta}}(\bar{z},e,\mu)\!=\!M c_{1}e\!-\!\sum_{i=1}^{r}c_{i}\bar{g}_{i}(\bar{z},e,v,w)$, $\bar{f}_{\bar{\xi}}(\bar{z},e,\tilde{\chi},\mu)=\lambda(b^{-1}(w)-1)e+G(\bar{z},e,v,w)+B\Psi\tilde{\eta}-\lambda\tilde{e}$, $h_{1}(X_{0},e,\mu)=b(w)\bar{\xi}_{2}+b(w)G_{1}(\bar{z},e,v,w)$,
$h_{i}(X_{i-1},\check{\xi}_{i},\tilde{\chi},\mu)=\lambda_{i}(e-b^{-1}(w)e+\bar{\xi}_{1})+\lambda_{i}\tilde{e}-\frac{\partial \vartheta_{i-1}(\check{\xi}_{i-1})}{\partial \check{\xi}_{i-1}}\dot{\check{\xi}}_{i-1}$ for $i=2,\cdots,r$.
To facilitate our analysis, we define the following systems:
\begin{equation}\label{system6}
\begin{split}
 & \dot{X}_{0} = F_{0}(X_{0},e,\tilde{\chi},\mu)\\
 & \dot{X}_{i} = F_{i}(X_{i},\check{\xi}_{i+1},\tilde{\chi},\mu),~i=1,\cdots,r-1\\
 & \dot{X}_{r} = F_{r}(X_{r},\tilde{\vartheta}_{r},\tilde{\chi},\mu)\\
\end{split}
\end{equation}
where $\check{\xi}_{r+1}=\tilde{\vartheta}_{r}$, $F_{0}(X_{0},e,\tilde{\chi},\mu)=\mbox{col}(\bar{f} (\bar{z} ,e,\mu),M \bar{\eta} +\bar{f}_{\bar{\eta}}(\bar{z},e,\mu),A_{o}\bar{\xi}+B\Psi\bar{\eta}+\bar{f}_{\bar{\xi}}(\bar{z},e,\tilde{\chi},\mu))$, $F_{1}(X_{1},\check{\xi}_{2},\tilde{\chi},\mu)=\mbox{col}(F_{0}(X_{0},e,\tilde{\chi},\mu),b(w)(\check{\xi}_{2}+\vartheta_{1})+h_{1}(X_{0},e,\mu))$,
$F_{i}(X_{i},\check{\xi}_{i+1},\tilde{\chi},\mu)=$\\
$\mbox{col}(F_{i-1}(X_{i-1},\check{\xi}_{i},\tilde{\chi},\mu),\check{\xi}_{i+1}+\vartheta_{i}+h_{i}(X_{i-1},\check{\xi}_{i},\tilde{\chi},\mu))$ for $i=2,\cdots,r-1$ and $F_{r}(X_{r},\tilde{\vartheta}_{r},\tilde{\chi},\mu)=\mbox{col}(F_{r-1}(X_{r-1},\check{\xi}_{r},\tilde{\chi},\mu),\vartheta_{r}+\tilde{\vartheta}_{r}+h_{r}(X_{r-1},\check{\xi}_{r},\tilde{\chi},\mu))$.
Clearly, the system $\dot{X_{r}} = F_{r}(X_{r},\tilde{\vartheta}_{r},\tilde{\chi},\mu)$ is equivalent to the system \eqref{system5}.

 \begin{Lemma}\label{Lemma1}
 Under Assumptions \ref{Ass2.1}-\ref{Ass2.4} and \ref{Ass3.1}, (i) there exists  a $\mathcal{C}^{1}$ function $U_{r}(X_{r})$ and some smooth positive functions $\rho_{i}(\cdot)$ with $i=1,\cdots,r$, such that, with
  \begin{equation}\label{varthetai1}
\begin{split}
\vartheta_{i}(\check{\xi}_{i})=-\rho_{i}(\check{\xi}_{i})\check{\xi}_{i},~i=1,\cdots,r,
\end{split}
\end{equation}
 for  any $\mu\in\Omega$, and any $X_{r}$, we have
 \begin{equation}\label{UXr1}
\begin{split}
\underline{\alpha}_{r}(\|X_{r}\|)\leq
U_{r}(X_{r})\leq\bar{\alpha}_{r}(\|X_{r}\|)
\end{split}
\end{equation}
\begin{equation}\label{dotUXr1}
\begin{split}
\frac{\partial U_{r}(X_{r})}{\partial X_{r}}F_{r}(X_{r},\tilde{\vartheta}_{r},\tilde{\chi},\mu)\leq\!-\alpha_{r}(\|X_{r}\|)\!+\!\tilde{\vartheta}_{r}^{2}+\pi_{r}(\tilde{\chi})
\end{split}
\end{equation}
where $\underline{\alpha}_{r}(\cdot)$,  $\bar{\alpha}_{r}(\cdot)$ and $\alpha_{r}(\cdot)$ are some class $\mathcal{K}_{\infty}$ functions with $\alpha_{r}(\|X_{r}\|)=-2\|X_{r}\|^{2}-\sigma^{2}\rho_{r}(\check{\xi}_{r})\check{\xi}_{r}^{2}$ for some real number $0<\sigma<1$, and $\pi_{r}(\cdot)$ is some positive definite function.
(ii) Consider the following event-triggered mechanism
\begin{equation}\label{trigger2}
\begin{split}
t_{k+1}\!=\!\inf\{t>t_{k}~|~&\tilde{\vartheta}_{r}^{2}(t)+\pi_{r}(\tilde{\chi}(t))-\sigma^{2}\rho_{r}(\check{\xi}_{r}(t))\check{\xi}_{r}^{2}(t)\geq\delta^{2}\}\\
\end{split}
\end{equation}
 where   $\delta$ is some positive real number. Then,  for  any $\mu\in\Omega$, and any $X_{r}$, we further have \begin{equation}\label{dotUXr4}
\begin{split}
\frac{\partial U_{r}(X_{r})}{\partial X_{r}}F_{r}(X_{r},\tilde{\vartheta}_{r},\tilde{\chi},\mu)\leq -\|X_{r}\|^{2},~ \forall~\|X_{r}\|\geq\delta.\\
\end{split}
\end{equation}
\end{Lemma}


\begin{Proof}
\textbf{Proof of Part (i):} The proof of this part consists of the following $r$ steps.

\textbf{Step} 1: From \eqref{system5}, we know that  the $X_{1}=\mbox{col}(\bar{z},\bar{\eta},\bar{\xi},e)$ subsystem is in the same form of equation (22) in \cite{PingHuang2013} with
$\zeta_{1}=\bar{z},\zeta_{2}=\mbox{col}(\bar{\eta},\bar{\xi})$, $x=e$, $y=\tilde{\chi}$, $\varphi_{1}=\bar{f} (\bar{z} ,e ,\mu)$, $\varphi_{2}=\mbox{col}(\bar{f}_{\bar{\eta}}(\bar{z},e,\mu),\bar{f}_{\bar{\xi}}(\bar{z},e,\tilde{\chi},\mu))$, $A=\left[
      \begin{array}{cc}
        M & 0 \\
        B\Psi & A_{o} \\
      \end{array}
    \right]$, $\phi=h_{1}(X_{0},e,\mu)$, $u=\check{\xi}_{2}+\vartheta_{1}(\check{\xi}_{1}),~\nu=\check{\xi}_{2}$.
Note that the matrix $\left[
      \begin{array}{cc}
        M & 0 \\
        B\Psi & A_{o} \\
      \end{array}
    \right]$ is Hurwitz, since the matrices $M$ and  $A_{o}$ are both Hurwitz.   Then, under Assumption \ref{Ass3.1}, by  Lemma 3.1 of \cite{PingHuang2013},  there exists a smooth positive function $\rho_{1}(\cdot)$ and a $\mathcal{C}^{1}$ function $U_{1}(X_{1})$ such that, for any $\mu\in\Omega$, and any $X_{1},\check{\xi}_{2},\tilde{\chi}$, with $\vartheta_{1}(\check{\xi}_{1})=-\rho_{1}(\check{\xi}_{1})\check{\xi}_{1}$,
 \begin{equation}\label{UX1}
\begin{split}
\underline{\alpha}_{1}(\|X_{1}\|)\leq
U_{1}(X_{1})\leq\bar{\alpha}_{1}(\|X_{1}\|)
\end{split}
\end{equation}
\begin{equation}\label{dotUX1}
\begin{split}
\frac{\partial U_{1}(X_{1})}{\partial X_{1}}F_{1}(X_{1},\check{\xi}_{2},\tilde{\chi},\mu)\leq -\|X_{1}\|^{2}+|\check{\xi}_{2}|^{2}+\pi_{1}(\tilde{\chi})\\
\end{split}
\end{equation}
where $\underline{\alpha}_{1}(\cdot)$ and $\bar{\alpha}_{1}(\cdot)$ are two class $\mathcal{K}_{\infty}$ functions,  and  $\pi_{1}(\cdot)$ is a smooth  positive definite function.

\textbf{Step $i$ ($2\leq i\leq r-1$)}: From \eqref{system5}, we also know that the subsystem $X_{i}=\mbox{col}(X_{i-1},\check{\xi}_{i})$ is in the same form of equation \eqref{zetasystem1} with
   \begin{equation*}
\begin{split}
 &\zeta_{1}=X_{i-1},~ \zeta_{2}=\check{\xi}_{i},~ \psi=\tilde{\chi},~b(\mu)=1\\
 &\varphi_{1}=F_{i-1}(X_{i-1},\check{\xi}_{i},\tilde{\chi},\mu),~\varphi_{2}=h_{i}(X_{i-1},\check{\xi}_{i},\tilde{\chi},\mu)\\
 &u=\check{\xi}_{i+1}+\vartheta_{i}(\check{\xi}_{i}),~\nu=\check{\xi}_{i+1}.\\
\end{split}
\end{equation*}
Assume that there exists a $\mathcal{C}^{1}$ function $U_{i-1}(X_{i-1})$ such that, for any $\mu\in\Omega$, and any $X_{i-1},\check{\xi}_{i},\tilde{\chi}$,
 \begin{equation}\label{UXi-1}
\begin{split}
\underline{\alpha}_{i-1}(\|X_{i-1}\|)\leq
U_{i-1}(X_{i-1})\leq\bar{\alpha}_{i-1}(\|X_{i-1}\|)
\end{split}
\end{equation}
\begin{equation}\label{dotUXi-1}
\begin{split}
&\frac{\partial U_{i-1}(X_{i-1})}{\partial X_{i-1}}F_{i-1}(X_{i-1},\check{\xi}_{i},\tilde{\chi},\mu)
\leq-\alpha_{i-1}(\|X_{i-1}\|)+|\check{\xi}_{i}|^{2}+\pi_{i-1}(\tilde{\chi})\\
\end{split}
\end{equation}
where $\pi_{i-1}(\cdot)$ is a smooth  positive definite function, $\underline{\alpha}_{i-1}(\cdot)$, $\bar{\alpha}_{i-1}(\cdot)$ and $\alpha_{i-1}(\cdot)$ are some class $\mathcal{K}_{\infty}$ functions  with  $\alpha_{i-1}(\cdot)$ satisfying $\lim_{s\rightarrow0^{+}}\sup(s^{2}/\alpha_{i-1}(s))<\infty$.
Then, by Lemma \ref{LemmaA1},  there exists a smooth positive function $\rho_{i}(\cdot)$ and a $\mathcal{C}^{1}$ function $U_{i}(X_{i})$ such that, for any $\mu\in\Omega$, and  any $X_{i},\check{\xi}_{i+1},\tilde{\chi}$,  with $\vartheta_{i}(\check{\xi}_{i}(t))=-\rho_{i}(\check{\xi}_{i}(t))\check{\xi}_{i}(t)$,
 \begin{equation}\label{UXi1}
\begin{split}
\underline{\alpha}_{i}(\|X_{i}\|)\leq
U_{i}(X_{i})\leq\bar{\alpha}_{i}(\|X_{i}\|)
\end{split}
\end{equation}
\begin{equation}\label{dotUXi1}
\begin{split}
&\frac{\partial U_{i}(X_{i})}{\partial X_{i}}F_{i}(X_{i},\check{\xi}_{i+1},\tilde{\chi},\mu)
\leq-\alpha_{i}(\|X_{i}\|)+|\check{\xi}_{i+1}|^{2}+\pi_{i}(\tilde{\chi})\\
\end{split}
\end{equation}
where $\pi_{i}(\cdot)$ is a smooth  positive definite function, $\underline{\alpha}_{i}(\cdot)$, $\bar{\alpha}_{i}(\cdot)$ and $\alpha_{i}(\cdot)$ are some class $\mathcal{K}_{\infty}$ functions with  $\alpha_{i}(\cdot)$ satisfying $\lim_{s\rightarrow0^{+}}\sup(s^{2}/\alpha_{i}(s))<\infty$.

\textbf{Step $r$}:
From step $r-1$, we know that there exists a $\mathcal{C}^{1}$ function $U_{r-1}(X_{r-1})$ such that \eqref{UXi1} and \eqref{dotUXi1} are satisfied with $i=r-1$. By the changing supply pair technique in \cite{Sontag1}, given any smooth function $\Delta_{r-1}(X_{r-1})\geq0$, there exists a $\mathcal{C}^{1}$ function $\bar{U}_{r-1}(X_{r-1})$, such that,
for any $\mu\in\Omega$, and any $X_{r-1},\check{\xi}_{r},\tilde{\chi}$,
 \begin{equation}\label{bUXr-1}
\begin{split}
\underline{\beta}_{r-1}(\|X_{r-1}\|)\leq \bar{U}_{r-1}(X_{r-1})\leq\bar{\beta}_{r-1}(\|X_{r-1}\|)
\end{split}
\end{equation}
 \begin{equation}\label{dotbUXr-1}
\begin{split}
&\frac{\partial \bar{U}_{r-1}(X_{r-1})}{\partial X_{r-1}}F_{r-1}(X_{r-1},\check{\xi}_{r},\tilde{\chi},\mu)
\leq-\Delta_{r-1}(X_{r-1})\|X_{r-1}\|^{2}\!\!+\!\bar{\phi}_{r-1}(\check{\xi}_{r})\check{\xi}_{r}^{2}\!\!+\!\bar{\pi}_{r-1}(\tilde{\chi})\|\tilde{\chi}\|^{2}\\
\end{split}
\end{equation}
where $\underline{\beta}_{r-1}(\cdot)$ and $\bar{\beta}_{r-1}(\cdot)$  are some known class $\mathcal{K}_{\infty}$ functions, $\bar{\phi}_{r-1}(\cdot)$ and $\bar{\pi}_{r-1}(\cdot)$ are some known smooth positive functions.

Note that $h_{r}(X_{r-1},\check{\xi}_{r},\tilde{\chi},\mu)$ is smooth and satisfies $h_{r}(0,0,0,\mu)=0$ for any $\mu\in\mathbb{R}^{n_{\mu}}$. Then, by applying Lemma 7.8 of \cite{Huang1}, there exist three smooth functions
$l_{r}(X_{r-1})$, $\kappa_{r}(\check{\xi}_{r})$ and $\omega_{r}(\tilde{\chi})$ satisfying $l_{r}(0)=0$, $\kappa_{r}(0)=0$ and $\omega_{r}(0)=0$  such that,  for  all $X_{r-1}$, $\check{\xi}_{r}$,  $\tilde{\chi}$ and all $\mu\in\Omega$,  $|h_{r}(X_{r-1},\check{\xi}_{r},\tilde{\chi},\mu)|\leq l_{r}(X_{r-1})+\kappa_{r}(\check{\xi}_{r})+\omega_{r}(\tilde{\chi})$,
which further implies $|h_{r}(X_{r-1},\check{\xi}_{r},\tilde{\chi},\mu)|^{2}=(l_{r}(X_{r-1})+\kappa_{r}(\check{\xi}_{r})+\omega_{r}(\tilde{\chi}))^{2}
\leq 3|l_{r}(X_{r-1})|^{2}+3|\kappa_{r}(\check{\xi}_{r})|^{2}+3|\omega_{r}(\tilde{\chi})|^{2}$.
Since $l_{r}(0)=0$, $\kappa_{r}(0)=0$ and $\omega_{r}(0)=0$, there exist some smooth positive functions  $\bar{l}_{r}(X_{r-1})$, $\bar{\kappa}_{r}(\check{\xi}_{r})$ and $\bar{\omega}_{r}(\tilde{\chi})$ such that, for all $X_{r-1}$,  $\check{\xi}_{r}$ and $\tilde{\chi}$, $3|l_{r}(X_{r-1})|^{2}\leq\bar{l}_{r}(X_{r-1})\|X_{r-1}\|^{2}$, $3|\kappa_{r}(\check{\xi}_{r})|^{2}\leq\bar{\kappa}_{r}(\check{\xi}_{r})\check{\xi}_{r}^{2}$ and $3|\omega_{r}(\tilde{\chi})|^{2}\leq\bar{\omega}_{r}(\tilde{\chi})\|\tilde{\chi}\|^{2}$.
As a result, for  all $X_{r-1}$, $\check{\xi}_{r}$,  $\tilde{\chi}$, and all $\mu\in\Omega$, we have
 \begin{equation}\label{hr4}
\begin{split}
&|h_{r}(X_{r-1},\check{\xi}_{r},\tilde{\chi},\mu)|^{2}
\leq\bar{l}_{r}(X_{r-1})\|X_{r-1}\|^{2}+\bar{\kappa}_{r}(\check{\xi}_{r})\check{\xi}_{r}^{2}+\bar{\omega}_{r}(\tilde{\chi})\|\tilde{\chi}\|^{2}.\\
\end{split}
\end{equation}

Let $V_{r}(\check{\xi}_{r})=\frac{1}{2}\check{\xi}_{r}^{2}$ and $\vartheta_{r}(\check{\xi}_{r})=-\rho_{r}(\check{\xi}_{r})\check{\xi}_{r}$ with $\rho_{r}(\cdot)$ being a smooth positive function, and let $\underline{\beta}_{r}(s)=\bar{\beta}_{r}(s)=\frac{1}{2}s^{2}$. Then $\underline{\beta}_{r}(|\check{\xi}_{r}|)\leq V_{r}(\check{\xi}_{r})\leq\bar{\beta}_{r}(|\check{\xi}_{r}|)$ for all $\check{\xi}_{r}$.
 Also, according to  \eqref{system5} and \eqref{hr4},  for all $\mu\in\Omega$, all $X_{r-1}$, all $\check{\xi}_{r}$ and all $\tilde{\chi}$,  we have
 \begin{equation}\label{Vr1}
\begin{split}
&\frac{\partial V_{r}(\check{\xi}_{r})}{\partial \check{\xi}_{r}}(\vartheta_{r}+\tilde{\vartheta}_{r}+h_{r}(X_{r-1},\check{\xi}_{r},\tilde{\chi},\mu))\\
=& \check{\xi}_{r}(\vartheta_{r}+\tilde{\vartheta}_{r}+h_{r}(X_{r-1},\check{\xi}_{r},\tilde{\chi},\mu))\\
\leq&-\!\rho_{r}(\check{\xi}_{r})\check{\xi}_{r}^{2}\!+\!\frac{1}{4}\check{\xi}_{r}^{2}\!+\!\tilde{\vartheta}_{r}^{2}\!+\!\frac{1}{4}\check{\xi}_{r}^{2}\!\!+\!\bar{l}_{r}(X_{r-1})\|X_{r-1}\|^{2}
+\bar{\kappa}_{r}(\check{\xi}_{r})\check{\xi}_{r}^{2}+\bar{\omega}_{r}(\tilde{\chi})\|\tilde{\chi}\|^{2}\\
=&-\!(\rho_{r}(\check{\xi}_{r})\!-\!\frac{1}{2}\!-\!\bar{\kappa}_{r}(\check{\xi}_{r}))\check{\xi}_{r}^{2}\!+\!\tilde{\vartheta}_{r}^{2}\!+\!\bar{l}_{r}(X_{r-1})\|X_{r-1}\|^{2}
+\bar{\omega}_{r}(\tilde{\chi})\|\tilde{\chi}\|^{2}.\\
\end{split}
\end{equation}

We further let $U_{r}(X_{r})=\bar{U}_{r-1}(X_{r-1})+V_{r}(\check{\xi}_{r})$. Then, by Lemma 11.3 of \cite{ChenZHuang2015}, we can choose some  class $\mathcal{K}_{\infty}$ functions $\underline{\alpha}_{r}(s)\leq\min\{\underline{\beta}_{r-1}(s/\sqrt{2}),\underline{\beta}_{r}(s/\sqrt{2})\}$ and $\bar{\alpha}_{r}(s)\geq\bar{\beta}_{r-1}(s)+\bar{\beta}_{r}(s)$ such that \eqref{UXr1} is satisfied. Also, according to  \eqref{dotbUXr-1} and \eqref{Vr1}, for all $\mu\in\Omega$, and all $X_{r-1}$,  $\check{\xi}_{r}$, $\tilde{\chi}$,  we have
\begin{equation}\label{dotUXr2}
\begin{split}
&\frac{\partial U_{r}(X_{r})}{\partial X_{r}}F_{r}(X_{r},\tilde{\vartheta}_{r},\tilde{\chi},\mu)  \\
\leq&\!-\!\Delta_{r-1}(X_{r-1})\|X_{r-1}\|^{2}\!\!+\!\bar{\phi}_{r-1}(\check{\xi}_{r})\check{\xi}_{r}^{2}\!\!+\!\bar{\pi}_{r-1}(\tilde{\chi})\|\tilde{\chi}\|^{2}
-(\rho_{r}(\check{\xi}_{r})\!-\!\frac{1}{2}\!-\!\bar{\kappa}_{r}(\check{\xi}_{r}))\check{\xi}_{r}^{2}\\
&+\tilde{\vartheta}_{r}^{2}\!+\!\bar{l}_{r}(X_{r-1})\|X_{r-1}\|^{2}+\bar{\omega}_{r}(\tilde{\chi})\|\tilde{\chi}\|^{2}\\
=&-(\Delta_{r-1}(X_{r-1})-\bar{l}_{r}(X_{r-1}))\|X_{r-1}\|^{2}
-(\rho_{r}(\check{\xi}_{r})\!-\!\frac{1}{2}\!-\!\bar{\kappa}_{r}(\check{\xi}_{r})-\bar{\phi}_{r-1}(\check{\xi}_{r}))\check{\xi}_{r}^{2}\\
&+\tilde{\vartheta}_{r}^{2}+(\bar{\pi}_{r-1}(\tilde{\chi})+\bar{\omega}_{r}(\tilde{\chi}))\|\tilde{\chi}\|^{2}.\\
\end{split}
\end{equation}
Choose  $\rho_{r}(\check{\xi}_{r})\geq \frac{1}{(1-\sigma^2)}(\frac{5}{2}+\bar{\kappa}_{r}(\check{\xi}_{r})+\bar{\phi}_{r-1}(\check{\xi}_{r}))$ and $\Delta_{r-1}(X_{r-1})\geq\bar{l}_{r}(X_{r-1}) +2$, and let $\pi_{r}(\tilde{\chi})=(\bar{\pi}_{r-1}(\tilde{\chi})+\bar{\omega}_{r}(\tilde{\chi}))\|\tilde{\chi}\|^{2}$.
Then we have
\begin{equation}\label{dotUXr3}
\begin{split}
&\frac{\partial U_{r}(X_{r})}{\partial X_{r}}F_{r}(X_{r},\tilde{\vartheta}_{r},\tilde{\chi},\mu)  \\
\leq&-2\|X_{r-1}\|^{2}-2\check{\xi}_{r}^{2}-\sigma^{2}\rho_{r}(\check{\xi}_{r})\check{\xi}_{r}^{2}+\tilde{\vartheta}_{r}^{2}+\pi_{r}(\tilde{\chi})\\
=&-2\|X_{r}\|^{2}-\sigma^{2}\rho_{r}(\check{\xi}_{r})\check{\xi}_{r}^{2}+\tilde{\vartheta}_{r}^{2}+\pi_{r}(\tilde{\chi})\\
\end{split}
\end{equation}
which implies (\ref{dotUXr1}).

\textbf{Proof of Part (ii)}: Note that, under the event-triggered mechanism \eqref{trigger2}, for any $t\in[t_{k},t_{k+1})$ with $k\in\mathbb{S}$, we have
\begin{equation}\label{tildevarthetar1}
\begin{split}
\tilde{\vartheta}_{r}^{2}(t)+\pi_{r}(\tilde{\chi}(t))\leq\sigma^{2}\rho_{r}(\check{\xi}_{r}(t))\check{\xi}_{r}^{2}(t)+\delta^{2}.
\end{split}
\end{equation}
Then, combining \eqref{dotUXr3} and \eqref{tildevarthetar1}, we have
\begin{equation}\label{dotUXr5}
\begin{split}
\frac{\partial U_{r}(X_{r})}{\partial X_{r}}F_{r}(X_{r},\tilde{\vartheta}_{r},\tilde{\chi},\mu)&\leq-2\|X_{r}\|^{2}+\delta^{2}\leq-\|X_{r}\|^{2},~\forall \|X_{r}\|\geq\delta.
\end{split}
\end{equation}
Thus the proof is completed.
\end{Proof}


\begin{Remark}\label{RemarkTrigger2}
It is of interest to explain the procedure for designing our triggering mechanism \eqref{trigger2}.
First, note that the term $\tilde{\vartheta}_{r}^{2}+\pi_{r}(\tilde{\chi})$ in the triggering mechanism \eqref{trigger2} is due to the fact that the system $\dot{X}_{r} = F_{r}(X_{r},\tilde{\vartheta}_{r},\tilde{\chi},\mu)$ is ISS with $\tilde{\vartheta}_{r}^{2}+\pi_{r}(\tilde{\chi})$ as the input.
Thus, it is naturally to require
      $\tilde{\vartheta}_{r}^{2}+\pi_{r}(\tilde{\chi}) \leq \delta$ for some $\delta >0$. Also, from the proof of Lemma \ref{Lemma1}, we know the functions $\tilde{\vartheta}_{r}$ and $\pi_{r}(\tilde{\chi})$ can be recursively constructed.
Second, the term  $\sigma^{2}\rho_{r}(\check{\xi}_{r})\check{\xi}_{r}^{2}$ is not necessary for the sake of the solvability of the problem.
  We use it for adjusting the transient performance and the triggering number during the initial stage, which is motivated by our previous work in \cite{LiuHuang2017b}. If we let $\sigma=0$, then the triggering number may become large during the initial stage.
Third, as will be discussed in Remark \ref{RemarkTheorem1a}, the term $\delta$ plays an important role in eliminating the Zeno behavior and also influences the steady-state tracking error and the triggering number during the steady stage.
\end{Remark}


\begin{Remark}
It is noted that not only the the event-triggered control law (\ref{baru2}) is more complicated than that in \cite{LiuHuang2017b}, but also
the event-triggered mechanism \eqref{trigger2} is more complicated than that in  \cite{LiuHuang2017b}.
This complexity is caused by two factors. First, unlike in \cite{LiuHuang2017b} where it suffices to stabilize the augmented system, here we need to
stabilize the extended augmented system \eqref{system4}, which contains some additional piecewise continuous functions $\tilde{\eta}$ and $\tilde{e}$.
Second, since the extended augmented system \eqref{system4} is of
 higher relative degree, the  functions
$\pi_{r}(\cdot)$ and $\vartheta_{r}(\cdot)$ have to be recursively constructed as detailed in the proof of
Lemma \ref{Lemma1}.
\end{Remark}

Proposition \ref{Proposition1} and Lemma \ref{Lemma1} lead to the following main result.

\begin{Theorem}\label{Theorem1}
 Under Assumptions \ref{Ass2.1}-\ref{Ass2.4} and \ref{Ass3.1}, for any $\epsilon>0$, there exists a $\delta>0$, such that 
  Problem \ref{Problem1} for the system \eqref{system1} is solvable by the following event-triggered output feedback control law
\begin{equation}\label{u3}
\begin{split}
 &u (t)=-\rho_{r}(\check{\xi}_{r}(t_{k}))\check{\xi}_{r}(t_{k})+\Psi \eta (t_{k}) \\
 &\dot{\hat{\xi}}(t)=A_{o}\hat{\xi}(t)+\lambda e(t_{k})+B(u(t)\!-\!\Psi \eta (t_{k}))\\
 &\dot{\eta} (t)=M \eta (t)+N u (t),~\forall t\in[t_{k},t_{k+1}),~k\in\mathbb{S}\\
 \end{split}
\end{equation}
under the output-based event-triggered mechanism \eqref{trigger2}.
\end{Theorem}
\begin{Proof}
Suppose that the solution $X_{r}(t)$ of the closed-loop system \eqref{system5} under the event-triggered mechanism \eqref{trigger2} is right maximally defined for all $t\in[0,T_{M})$ with $0<T_{M}\leq\infty$. Then, based on Lemma \ref{Lemma1}, we have
\begin{equation} \label{bounded}
\|X_{r} (t)\| \leq \max \{\delta, \underline{\alpha}_{r}^{-1}(\bar{\alpha}_{r}(\|X_{r}(0)\|))\}, \forall t\in[0,T_M).
\end{equation}

We first consider the case that $\mathbb{S}=\mathbb{N}$, i.e., there are infinite many triggering times. If we show $\lim_{k\rightarrow\infty}t_{k}=\infty$, then $T_{M}$ must be equal to $\infty$.
For this purpose,  note that, for any $t\in[t_{k},t_{k+1})$ with $k\in\mathbb{S}$,
 \begin{equation}\label{dottildevarthetai1}
\begin{split}
&\frac{d\big(\tilde{\vartheta}_{r}^{2}(t)+\pi_{r}(\tilde{\chi}(t))\big)}{d t}\\
=&2\tilde{\vartheta}_{r}(t)\dot{\tilde{\vartheta}}_{r}(t)+\frac{d \pi_{r}(\tilde{\chi})}{d \tilde{\chi}}\dot{\tilde{\chi}}(t)\\
=&-2\tilde{\vartheta}_{r}(t)\frac{d\vartheta_{r}(\check{\xi}_{r}) }{d \check{\xi}_{r}}\dot{\check{\xi}}_{r}(t)-\frac{d \pi_{r}(\tilde{\chi})}{d \tilde{\chi}}\dot{\chi}(t)\\
\end{split}
\end{equation}
where, for convenience, we still use $\frac{d(\tilde{\vartheta}_{r}^{2}(t)+\pi_{r}(\tilde{\chi}(t)))}{d t}$, $\dot{\tilde{\vartheta}}_{r}(t)$, $\dot{\tilde{\chi}}(t)$, $\dot{\check{\xi}}_{r}(t)$ and $\dot{\chi}(t)$ to denote their right derivatives at the triggering time instant $t_{k}$.
From \eqref{bounded}, we know $X_{r}(t)$  is bounded for all $t \in[0,T_{M})$. Together with the definitions of $\tilde{\vartheta}_{r}$ and $\tilde{\chi}$, we know $\tilde{\vartheta}_{r}$ and $\tilde{\chi}$ are also bounded for all $t \in[0,T_{M})$. Then, from $\dot{X}_{r} = F_{r}(X_{r},\tilde{\vartheta}_{r},\tilde{\chi},\mu)$, we know $\dot{X}_{r}(t)$ is also bounded for all $t \in[0,T_{M})$, which further implies that $\dot{\tilde{\vartheta}}_{r}(t)$, $\dot{\tilde{\chi}}(t)$, $\dot{\check{\xi}}_{r}(t)$ and $\dot{\chi}(t)$ are bounded for all $t \in[0,T_{M})$.
Then we conclude that there exists a positive constant $c_{0}$ depending on $\delta$ and $X_{r}(0)$ such that
$\frac{d\big(\tilde{\vartheta}_{r}^{2}(t)+\pi_{r}(\tilde{\chi}(t))\big)}{d t}\leq c_{0}$ for all $t\in[0,T_{M})$.

On the other hand,  according to the definitions of $\tilde{\chi}(t)$ and $\tilde{\vartheta}(t)$,  we know
$\tilde{\chi}(t_{k})=\chi (t_{k})-\chi (t_{k})=0$ and $\tilde{\vartheta}(t_{k})=\vartheta (\check{\xi}(t_{k}))-\vartheta (\check{\xi}(t_{k}))=0$ for all $k\in\mathbb{S}$,
and, from \eqref{trigger2}, we have
$\lim_{t\rightarrow t_{k+1}^{-}}(\tilde{\vartheta}_{r}^{2}(t)+\pi_{r}(\tilde{\chi}(t)))\geq\delta^{2}$ for all $k\in\mathbb{S}$.
As a result, we conclude that $\delta^{2}\leq c_{0}(t_{k+1}-t_{k})$ for all $k\in\mathbb{S}$,
which implies $t_{k+1}-t_{k}\geq \frac{\delta^{2}}{c_{0}}$ for all $k\in\mathbb{S}$.
Thus $\lim_{k\rightarrow\infty}t_{k}=\infty$ and the Zeno behavior does not happen. Then $T_{M}=\infty$. 

Next, we consider the case that $\mathbb{S}=\{1,2,\cdots,k^{*}\}$ with $k^{*}$ a positive integer.
For this case, the closed-loop system \eqref{system5} reduces to a time-invariant continuous-time system for all $t\geq t_{k^{*}}$. Then, according to \eqref{bounded} and  the definition of $T_{M}$, we have $T_{M}=\infty$.

Since the solution $X_{r}(t)$ of the closed-loop system \eqref{system5} exists for all $t\in[0,\infty)$,  by applying Theorem 4.18 of \cite{Khalil1} and Lemma \ref{Lemma1} here,  we have that $X_{r}(t)$ is globally ultimately bounded with the ultimate bound $ d(\delta)=\underline{\alpha}_{r}^{-1}(\bar{\alpha}_{r}(\delta))$, i.e., $\lim_{t\rightarrow \infty} \sup \|X_{r}(t)\|\leq d(\delta)=\underline{\alpha}_{r}^{-1}(\bar{\alpha}_{r}(\delta))$.
Note that  $d(\cdot)$ is  an invertible class $\mathcal{K}_{\infty}$ function since $\underline{\alpha}_{r}(\cdot)$ and $\bar{\alpha}_{r}(\cdot)$ are both invertible class $\mathcal{K}_{\infty}$ functions. For any $\epsilon>0$, let $\delta=d^{-1}(\epsilon)=\bar{\alpha}_{r}^{-1} (\underline{\alpha}_{r}( \epsilon))$. Then  we have $\lim_{t\rightarrow\infty}\sup\|X_{r}(t)\|\leq\epsilon$.
That is to say, the control law \eqref{baru2} together with the event-triggered mechanism \eqref{trigger2} solves the global robust practical stabilization problem for the extended augmented system \eqref{system4}.

Using Proposition \ref{Proposition1} completes the proof.
\end{Proof}

\begin{Remark}\label{RemarkTheorem1a}
As remarked in \cite{LiuHuang2017b},  letting $\delta = 0$ gives a special case of Theorem \ref{Theorem1} where $\lim_{t\rightarrow\infty}|e(t)| = \lim_{t\rightarrow\infty}\|X_{r}(t)\|=0$.
Nevertheless, in this case,  we cannot guarantee the prevention of the Zeno behavior. That is why we have introduced the positive constant $\delta$ in the triggering mechanism \eqref{trigger2}.
 In addition, from $\delta = \bar{\alpha}_{r}^{-1} (\underline{\alpha}_{r}( \epsilon))$, we can conclude that a smaller $\delta$ leads to a smaller steady-state tracking error. However, the price for a smaller $\delta$ is that the triggering number may become larger during the steady stage.
\end{Remark}

\begin{Remark}\label{RemarkTheorem1c}
The control law \eqref{u3} lends itself to the  following  digital implementation:
\begin{equation}\label{u4}
\begin{split}
 &u (t)=-\rho_{r}(\check{\xi}_{r}(t_{k}))\check{\xi}_{r}(t_{k})+\Psi \eta (t_{k}) \\
 &\hat{\xi}(t_{k+1})=\textbf{e}^{A_{o}(t_{k+1}-t_{k})}\hat{\xi}(t_{k})+\int_{t_{k}}^{t_{k+1}}\textbf{e}^{A_{o}(t_{k+1}-\tau)}d\tau(\lambda e(t_{k})-B\rho_{r}(\check{\xi}_{r}(t_{k}))\check{\xi}_{r}(t_{k}))\\
 &\eta(t_{k+1})=\textbf{e}^{M(t_{k+1}-t_{k})}\eta(t_{k})+\int_{t_{k}}^{t_{k+1}}\textbf{e}^{M(t_{k+1}-\tau)}d\tau N(-\rho_{r}(\check{\xi}_{r}(t_{k}))\check{\xi}_{r}(t_{k})+\Psi \eta (t_{k}))\\
 \end{split}
\end{equation}
for any $t\in[t_{k},t_{k+1})$ with $k\in\mathbb{S}$. 
 %
\end{Remark}


\section{An Example}\label{Example}
Consider the controlled hyper-chaotic Lorenz systems  \cite{XuHuang2010c} described as follows
\begin{equation}\label{Lorenzsystem1}
\begin{split}
 &\dot{z}_{1}=a_{1}z_{1}+a_{2}x_{1} \\
 &\dot{z}_{2}=a_{3}z_{2}+z_{1}x_{1} \\
 &\dot{x}_{1} =x_{2}+a_{4}z_{1}+a_{5}x_{1} -z_{1}z_{2} \\
&\dot{x}_{2}=bu+a_{6}z_{1}\\
 &y= x_{1}\\
 \end{split}
\end{equation}
where  $a \triangleq\mbox{col} (a_{1},\cdots,a_{6},b )$ is a constant parameter vector satisfying $a_{1}<0$, $a_{3}<0$ and $b >0$.
To account for the uncertain parameter vector, we let $a =\bar{a} +w$, where
$\bar{a} =\mbox{col}(\bar{a}_{1},\cdots,\bar{a}_{6},\bar{b} )=\mbox{col}(-8,1,-6,2,-1,-2,1)$
represents the nominal value of $a$, and $w =\mbox{col}(w_{1},\cdots,w_{7})$ represents the uncertainty of $a$. It is easy to see that the system \eqref{Lorenzsystem1} is in the form \eqref{system1} with relative degree $r=2$.

Let $v=\mbox{col}(v_{1},v_{2})\in\mathbb{R}^{2}$. Then we define an exosystem system in the form (\ref{exosystem1}) as follows
\begin{equation}\label{exoysytem2}
\begin{split}
\left[
     \begin{array}{c}
       \dot{v}_{1} \\
       \dot{v}_{2} \\
     \end{array}
   \right]=\left[
  \begin{array}{cc}
    0 & 1\\
    -1 & 0 \\
  \end{array}\right]\left[
     \begin{array}{c}
       v_{1} \\
       v_{2} \\
     \end{array}
   \right],~y_{0}=v_{1}.
 \end{split}
\end{equation}
 Clearly, Assumptions \ref{Ass2.1} is satisfied.
 We assume that $w\in\mathbb{W}=\{w\;|\;w\in\mathbb{R}^{7}, |w_{i}|\leq1,~i=1,\cdots,7\}$, and $v\in\mathbb{V}=\{v\;|\;v\in\mathbb{R}^{2}, |v _{i}|\leq1,~i=1,2\}$.

Define the tracking error $e=y-y_{0}$. Then, as shown in \cite{XuHuang2010c}, the regulator equations associated with \eqref{Lorenzsystem1} and \eqref{exoysytem2} are solvable and the solutions are given as follows
\begin{equation*}\label{regulator2}
\begin{split}
&\textbf{z}_{1}(v,w)=r_{11}(w)v_{1}+r_{12}(w)v_{2}\\
&\textbf{z}_{2}(v,w)=r_{21}(w)v_{1}^{2}+r_{22}(w)v_{1}v_{2}+r_{23}(w)v_{2}^{2}\\
&\textbf{x}_{1} (v,w)=v_{1}\\
&\textbf{x}_{2} (v,w)=r_{31}(w)v_{1}\!+\!r_{32}(w)v_{2}\!+\!r_{33}(w)v_{1}^{3}\!+\!r_{34}(w)v_{1}^{2}v_{2}+r_{35}(w)v_{1}v_{2}^{2}+r_{36}(w)v_{2}^{3}\\
&\textbf{u} (v,w)=r_{41}(w)v_{1}\!+\!r_{42}(w)v_{2}\!+\!r_{43}(w)v_{1}^{3}\!+\!r_{44}(w)v_{1}^{2}v_{2}+r_{45}(w)v_{1}v_{2}^{2}+r_{46}(w)v_{2}^{3}\\
 \end{split}
\end{equation*}
where the coefficients are the same as those given in \cite{XuHuang2010c}. As a result, Assumptions \ref{Ass2.3} and \ref{Ass2.4} are both satisfied.
It can be verified that
\begin{equation}
\begin{split}
\dfrac{d^{4}\textbf{u} (v,w)}{dt^{4}}+10\dfrac{d^{2}\textbf{u} (v,w)}{dt^{2}}+9\textbf{u} (v,w)=0.
 \end{split}
\end{equation}
Thus we have 
\begin{equation*}
\begin{split}
&\Phi =\left[
                                                                                                              \begin{array}{cccc}
0 & 1 & 0 & 0 \\
0 & 0 & 1 & 0 \\
0& 0 & 0 & 1 \\
-9& 0 & -10 & 0 \\
                                                                                                              \end{array}
                                                                                                            \right]
,\ \ \Gamma =\left[
                  \begin{array}{c}
                    1 \\
                    0 \\
                    0 \\
                    0 \\
                  \end{array}
                \right]^{T}.
 \end{split}
\end{equation*}
Then we can define an internal model of the form \eqref{doteta1} with
\begin{equation*}
\begin{split}
M =\left[
                                                                                                              \begin{array}{cccc}
0 & 1 & 0 & 0 \\
0 & 0 & 1 & 0 \\
0 & 0 & 0 & 1 \\
-4 & -12 & -13 & -6 \\
                                                                                                              \end{array}
                                                                                                            \right],\
                                                                                                            N =\left[
                  \begin{array}{c}
                    0 \\
                    0 \\
                    0 \\
                    1 \\
                  \end{array}
                \right].
\end{split}
\end{equation*}
Solving the Sylvester equation $T \Phi -M T =N \Gamma$ gives $\Psi =\Gamma T^{-1} =[-5, 12, 3, 6]$.

Perform the coordinate transformations \eqref{transformation1} and \eqref{transformation2}.  Then we obtain the following augmented system
 \begin{equation}\label{Lorenzsystem2}
\begin{split}
  \dot{\bar{z}} =&\bar{f} (\bar{z} ,e,v,w)\\
  \dot{\bar{\eta}} =&M \bar{\eta} +M c_{1}e-c_{1}\bar{g}_{1}(\bar{z},e,v,w)-c_{2}\bar{g}_{2}(\bar{z},e,v,w)\\
  \dot{\xi} =&A_{c}\xi+B\Psi\bar{\eta}+\bar{G}(\bar{z},e,v,w)+B\bar{u}+B\Psi\tilde{\eta}\\
\end{split}
\end{equation}
where
\begin{equation*}
\begin{split}
&A_{c}=\left[
         \begin{array}{cc}
           0 & 1 \\
           0 & 0 \\
         \end{array}
       \right],~B=\left[
        \begin{array}{c}
          0 \\
          1 \\
        \end{array}
      \right],~U_{d}=\left[
                \begin{array}{cc}
                  1 & 0 \\
                  -d_{1} & 1 \\
                \end{array}
              \right]\\
&c_{1}=b^{-1}MN,~c_{2}=b^{-1}N\\
&d_{1}=\Psi N,~d_{2}=\Psi MN,~\bar{z}=\mbox{col}(\bar{z}_{1},\bar{z}_{2})\\
&\bar{f} (\bar{z} ,e ,v,w)= \left[
                              \begin{array}{c}
                                a_{1}\bar{z}_{1}+a_{2}e \\
                                a_{3}\bar{z}_{2}+(\bar{z}_{1}+\textbf{z}_{1})(e+v_{1})-\textbf{z}_{1}v_{1} \\
                              \end{array}
                            \right]\\
&\bar{g}_{1}(\bar{z},e,v,w)=a_{4}\bar{z}_{1}+a_{5}e-\bar{z}_{1}\bar{z}_{2}-\bar{z}_{1}\mathbf{z}_{2}-\mathbf{z}_{1}\bar{z}_{2}\\
&\bar{g}_{2}(\bar{z},e,v,w)=a_{6}\bar{z}_{1}\\
&G_{1}(\bar{z},e,v,w)=b^{-1}(d_{1}e+\bar{g}_{1}(\bar{z},e,v,w))\\
&G_{2}(\bar{z},e,v,w)\!=\!b^{-1}\big(d_{2}e\!-\!d_{1}\bar{g}_{1}(\bar{z},e,v,w)\!+\!\bar{g}_{2}(\bar{z},e,v,w)\big)\\
&G(\bar{z},e,v,w)=\mbox{col}(G_{1}(\bar{z},e,v,w),G_{2}(\bar{z},e,v,w)).\\
\end{split}
\end{equation*}
Design an observer of the form \eqref{observer1} and define the observer error $\bar{\xi}=\xi-\hat{\xi}$. Attach the observer error dynamics to the augmented system \eqref{Lorenzsystem2} and replace the state variable $\xi$ by $\mbox{col}(e,\hat{\xi}_{2},\cdots,\hat{\xi}_{r})$. Then we get the following extended augmented system
\begin{equation}\label{Lorenzsystem3}
\begin{split}
  \dot{\bar{z}} =&\bar{f} (\bar{z} ,e,v,w)\\
  \dot{\bar{\eta}} =&M \bar{\eta} +M c_{1}e-c_{1}\bar{g}_{1}(\bar{z},e,v,w)-c_{2}\bar{g}_{2}(\bar{z},e,v,w)\\
  \dot{\bar{\xi}}=&A_{o}\bar{\xi}+\lambda(b^{-1}-1)e+B\Psi\bar{\eta}+G(\bar{z},e,v,w)+B\Psi\tilde{\eta}-\lambda\tilde{e}\\
\dot{e}=&b\bar{\xi}_{2}+b\hat{\xi}_{2}+bG_{1}(\bar{z},e,v,w)\\
\dot{\hat{\xi}}_{2}\!=&\bar{u}+\lambda_{2}(e-\hat{\xi}_{1})+\lambda_{2}\tilde{e}\\
\end{split}
\end{equation}
where
\begin{equation*}
\begin{split}
&\lambda\!=\!\left[
             \begin{array}{c}
               \lambda_{1} \\
               \lambda_{2} \\
             \end{array}
           \right]\!=\!\left[
             \begin{array}{c}
               2 \\
               2 \\
             \end{array}
           \right]\!,~A_{o}\!=\!\left[
            \begin{array}{cc}
              -2 & 1 \\
              -2 & 0 \\
            \end{array}
          \right].\\
\end{split}
\end{equation*}

In order to verify Assumption \ref{Ass3.1}, we choose the Lyapunov function candidate for the $\bar{z}$-subsystem  as follows:
\begin{equation}
\begin{split}
V_{0}(\bar{z} )=\frac{\hbar }{2}\bar{z}_{1}^{2}+\frac{\hbar }{4}\bar{z}_{1}^{4}+\frac{1}{2}\bar{z}_{2}^{2}
\end{split}
\end{equation}
for some sufficiently large $\hbar >0$. Since $V_{0}(\bar{z})$ is a $\mathcal{C}^{1}$ positive definite and radially unbounded function, there exist two class $\mathcal{K}_{\infty}$ functions $\underline{\alpha}_{0}(\cdot)$ and $\bar{\alpha}_{0}(\cdot)$ such that condition \eqref{Vz0} is satisfied.
 It is also possible to show that, for all $\mu\in\Omega$, all $\bar{z}$ and all $e$,
 \begin{equation}
\begin{split}
  \frac{\partial V_{0}(\bar{z})}{\partial \bar{z}} \bar{f} (\bar{z} ,e,v,w) \leq &-\ell_{1}\bar{z}_{1}^{2}-\ell_{2}\bar{z}_{1}^{4}-\ell_{3}\bar{z}_{2}^{2}+\ell_{4}e ^{2}+\ell_{5}e ^{4}
\end{split}
\end{equation}
for some positive constants $\ell_{1},\cdots,\ell_{5}$. Thus Assumption \ref{Ass3.1} is satisfied.

Then, by applying Theorem \ref{Theorem1} and following the detailed procedures in Section \ref{MR},   we can design the following event-triggered output feedback  control law
\begin{equation}\label{u5}
\begin{split}
 &u (t)=\vartheta_{2}(\check{\xi}_{2}(t_{k}))+\Psi \eta (t_{k}) \\
 &\dot{\hat{\xi}}(t)=A_{o}\hat{\xi}(t)+\lambda e(t_{k})+B(u(t)-\Psi \eta (t_{k}))\\
 &\dot{\eta} (t)=M \eta (t)+N u (t),~\forall t\in[t_{k},t_{k+1}),~k\in\mathbb{S}\\
 \end{split}
\end{equation}
and the following output-based event-triggered mechanism
\begin{equation}\label{trigger3}
\begin{split}
&t_{k+1}=\inf\{t>t_{k}~|~\tilde{f}(\tilde{\vartheta}_{2}(t),\tilde{\eta}(t),\tilde{e}(t),\tilde{\xi}_{1}(t))\geq \sigma^{2}|\vartheta_{2}(\check{\xi}_{2}(t))\check{\xi}_{2}(t)|+\delta^{2}\}\\
\end{split}
\end{equation}
where  $\check{\xi}_{2}(t)=\hat{\xi}_{2}(t)+6(e^{6}(t)+1)$, $\vartheta_{2}(\check{\xi}_{2}(t))=-12(\check{\xi}_{2}^{2}(t)+1)\check{\xi}_{2}(t)$, $\tilde{f}(\tilde{\vartheta}_{2}(t),\tilde{\eta}(t),\tilde{e}(t),\tilde{\xi}_{1}(t))=\tilde{\vartheta}_{2}^{2}(t)+5\|B\Psi\tilde{\eta}(t)-\lambda\tilde{e}(t)\|^{4}
+|\lambda_{2}\tilde{e}(t)|^{2}$, $\sigma=0.4$, and $\delta=0.1$ or $0.01$.



\begin{figure}[H]
\centering
\includegraphics[scale=0.7]{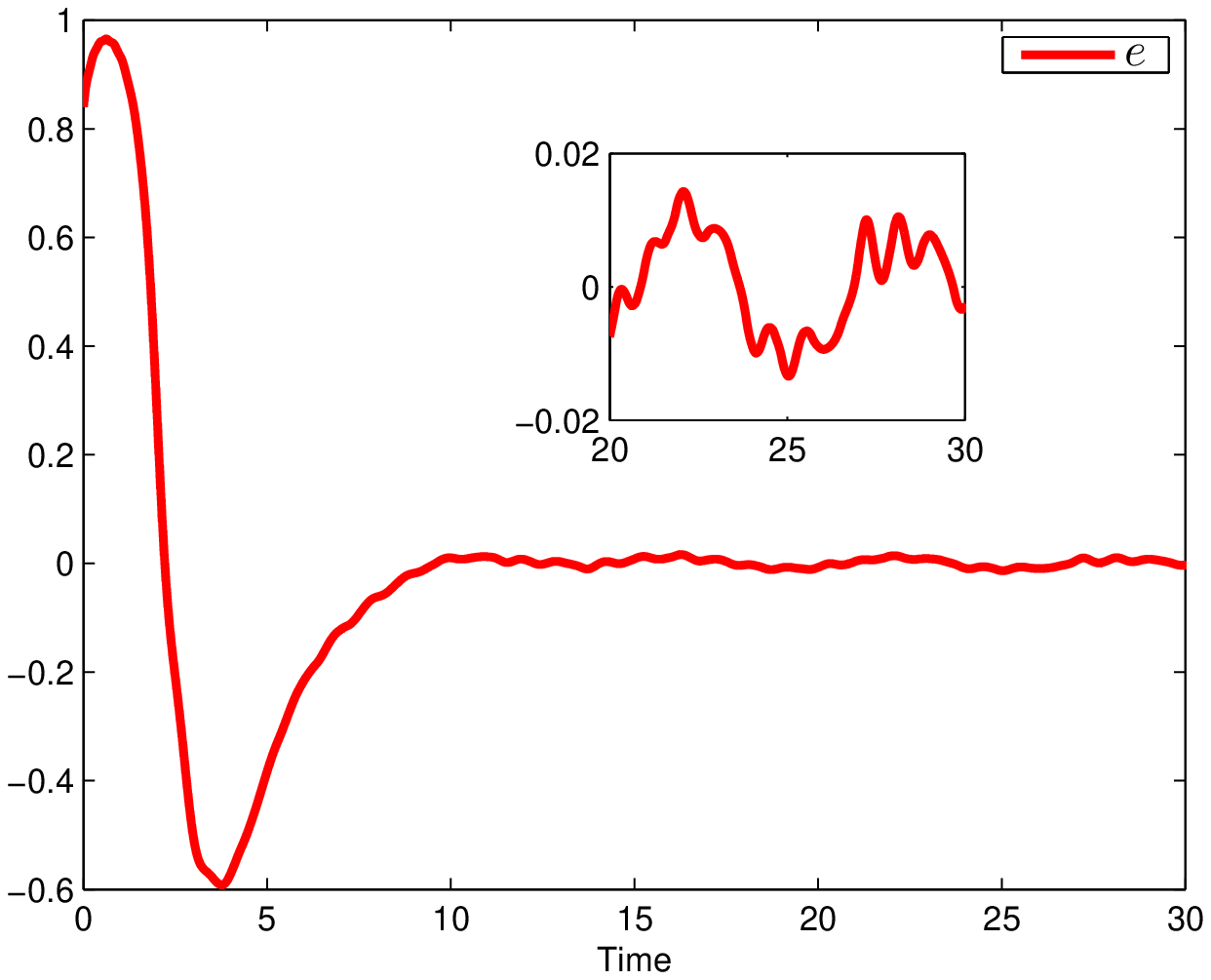}
\caption{Tracking error for $\delta=0.1$. } \label{error1}
\end{figure}
\begin{figure}[H]
\centering
\includegraphics[scale=0.7]{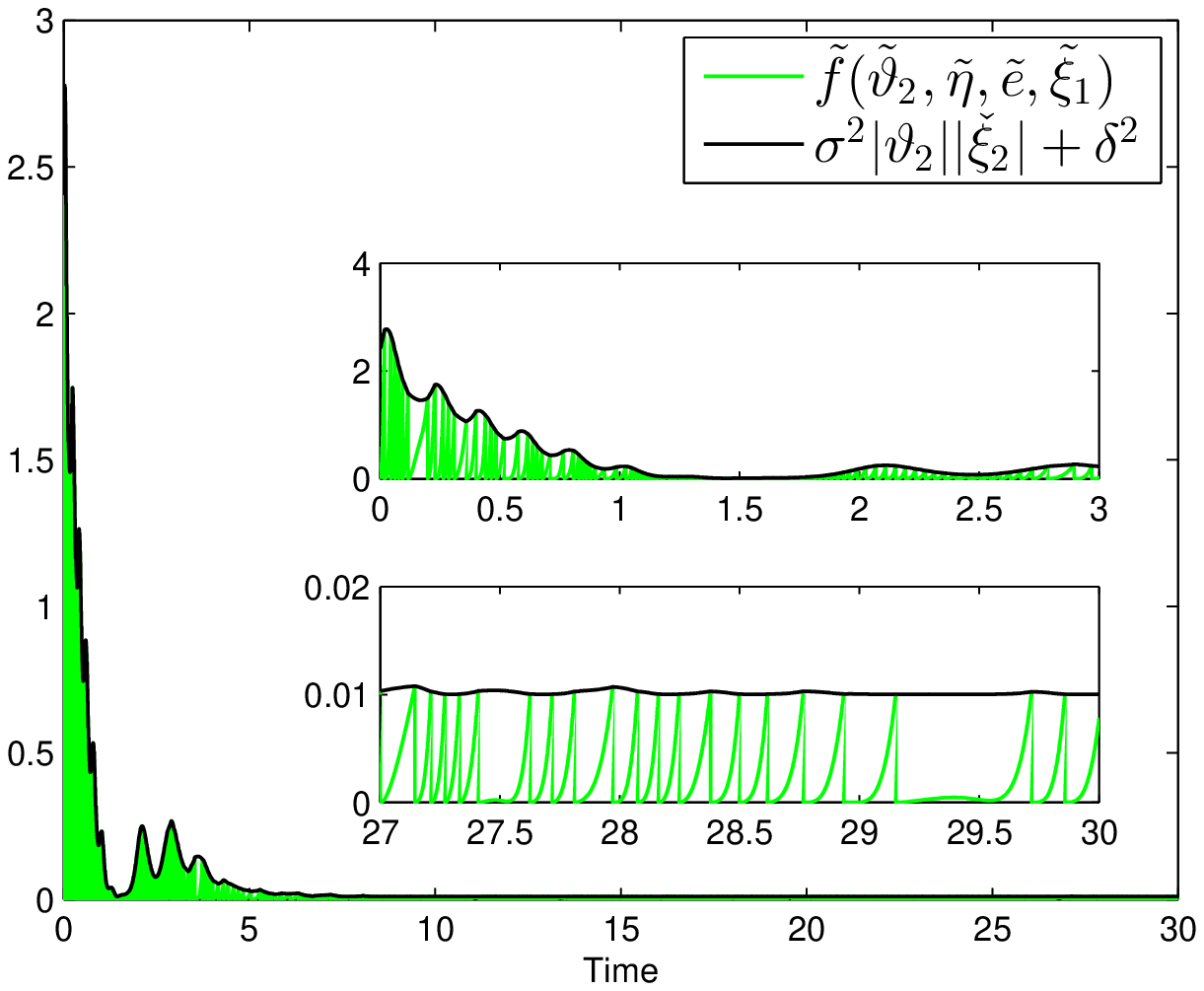}
\caption{Event-triggered  condition for $\delta=0.1$.} \label{condition1}
\end{figure}
\begin{figure}[H]
\centering
\includegraphics[scale=0.7]{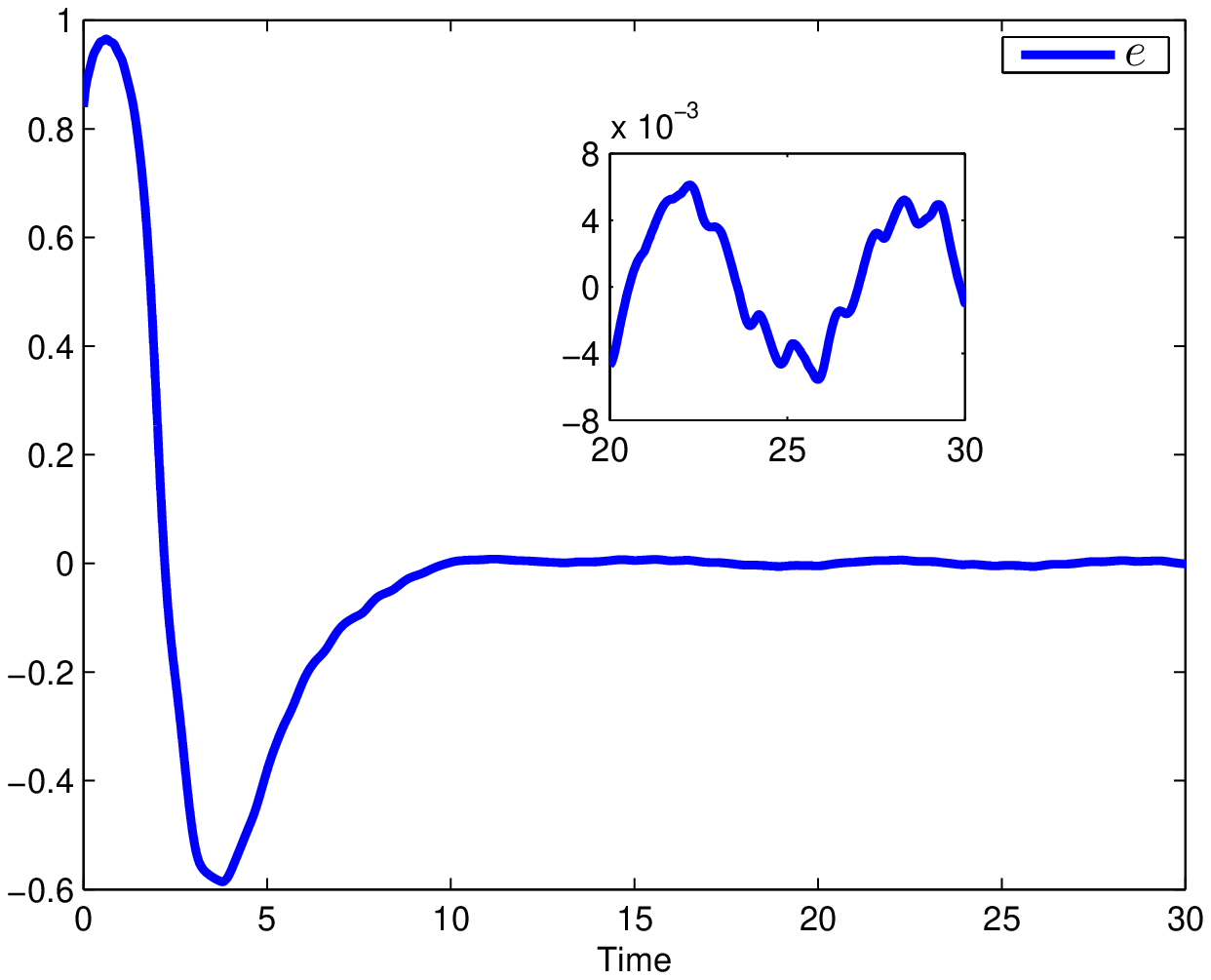}
\caption{Tracking error for $\delta=0.01$. } \label{error2}
\end{figure}
\begin{figure}[H]
\centering
\includegraphics[scale=0.7]{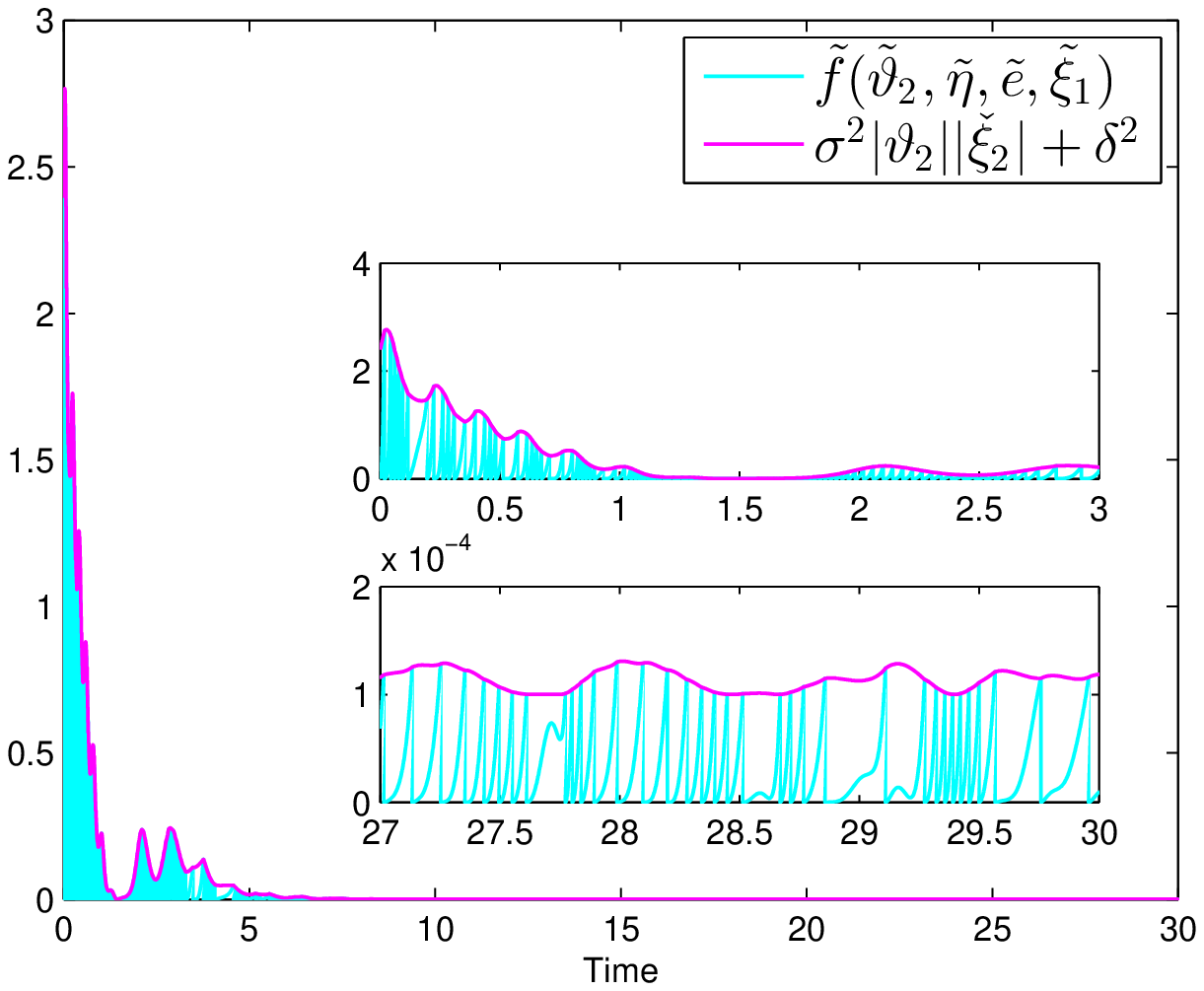}
\caption{Event-triggered  condition for $\delta=0.01$.} \label{condition2}
\end{figure}

\begin{table}[H]
  \begin{center}
    \caption{Event-triggered numbers.}\label{table}
  \scalebox{1.2}{
    \begin{tabular}{|c|c|c|} \hline
      \multirow{2}{*}{Design parameters} & \multirow{2}{*}{Time} & \multirow{2}{*}{Triggering numbers } \\
       &  &  \\ \hline
      \multirow{2}{*}{ $\sigma=0.4$, $\delta=0.1$} & \multirow{2}{*}{0--30s} & \multirow{2}{*}{271} \\
       &  &  \\ \hline
      \multirow{2}{*}{ $\sigma=0.4$, $\delta=0.01$} & \multirow{2}{*}{0--30s} & \multirow{2}{*}{478} \\
       &  &  \\ \hline
    \end{tabular}}
  \end{center}
\end{table}

Simulations are performed with  $w=[0.5,-0.4,0.1,-0.3$, $0.2,-0.3,0.4]^{T}$
 and the following initial conditions  
\begin{equation*}
\begin{split}
 &v(0)=[-0.34,-0.94]^{T},~z(0)=[0.13,   -0.67]^{T}\\
 &x(0)=[0.50,    0.30]^{T},~\hat{\xi}(0)=[ -1.40,   -5.96]^{T}\\
 &\eta(0)=[-0.35,   1.50,  -1.49,    0.31]^{T}.\\
 \end{split}
\end{equation*}

 The event-triggered numbers for both $\delta=0.1$ and $\delta=0.01$ are shown in Table \ref{table}. The tracking errors for $\delta=0.1$ and $\delta=0.01$ are shown in Figures \ref{error1} and \ref{error2}, respectively. The event-triggered conditions for $\delta=0.1$ and $\delta=0.01$ are shown in Figures \ref{condition1} and \ref{condition2}, respectively.
It can be seen from Table I that the triggering number for $\delta=0.01$ is greater than that for $\delta=0.1$.
Also, from Figures  \ref{error1} and \ref{error2}, we can see that $\lim_{t\rightarrow\infty}\sup |e (t)|\leq 0.02$ for $\delta=0.1$, and $\lim_{t\rightarrow\infty}\sup |e (t)|\leq 0.008$ for $\delta=0.01$.
Thus these simulation results illustrate that a larger $\delta$ leads to
less triggering number but leads to larger steady-state tracking error, which coincides with the viewpoint noted in Remark \ref{RemarkTheorem1a}.

\section{Conclusion}\label{Conclusion}
In this paper, we have studied  the event-triggered  global robust practical output regulation problem for nonlinear systems in output feedback form with any relative degree. The problem contains the problem in \cite{LiuHuang2017b} as a special case.
The higher relative degree of the plant incurs at least three specific challenges. First, since the partial state ${\xi}$ of the augmented system \eqref{system3} is not available for control,  we need to design an extra observer relying on the sampled error output $e(t_{k})$  to estimate the partial state ${\xi}$ of \eqref{system3}, thus entailing the stabilization of a much more complicated  extended augmented system. Second, since the  extended augmented system \eqref{system4} contains some aperiodically sampled quantities ${\eta} (t_{k})$ and ${e} (t_{k})$, $k\in\mathbb{S}$,
we need to establish the technical Lemma \ref{LemmaA1} to recursively design the dynamic output feedback control law and the output-based event-triggered mechanism. Third, due to the complexity of the closed-loop system, the stability analysis is also more complex than the unity relative degree case in \cite{LiuHuang2017b}.
It would be interesting to further
consider the event-triggered cooperative global robust practical output regulation problem for nonlinear multi-agent systems in output feedback form.


\begin{appendix}
\subsection{Proof of Proposition \ref{Proposition1}}\label{appendix1}
\begin{Proof}
First, note that $|e(t)|\leq\|\bar{x}_{c}(t)\|$. Together with \eqref{pcgs1}, we have
\begin{equation}\label{pcgs2}
\lim_{t \to \infty}\sup|e(t)|\leq \epsilon.
\end{equation}
which means that Property 2) of Problem \ref{Problem1} is satisfied.

 Next, we only need to show that Property 1) of Problem \ref{Problem1} is also satisfied.  For this purpose, we denote the state  of the closed-loop system composed of  \eqref{system1} and \eqref{u2} under the triggering mechanism \eqref{trigger1} by $x_{c}(t)=\mbox{col}(z(t),x(t),\eta(t),\hat{\xi}(t))$.

Since $\bar{x}_{c}(t)$ exists and is bounded for all $t\in[0,\infty)$, we know that $\bar{z}(t),\bar{\eta}(t),\bar{\xi}(t),e(t),\hat{\xi}_{2}(t)$, $\cdots,\hat{\xi}_{r}(t)$ are all bounded for all $t\in[0,\infty)$. Besides, due to the fact that $\hat{\xi}_{1}=\xi_{1}-\bar{\xi}_{1}=b^{-1}(w)e-\bar{\xi}_{1}$, we have that $\hat{\xi}_{1}(t)$ is bounded for all $t\in[0,\infty)$. Thus $\hat{\xi}(t)=\mbox{col}(\hat{\xi}_{1}(t),\cdots,\hat{\xi}_{r}(t))$ is bounded for all $t\in[0,\infty)$.
  According to the coordinate transformations \eqref{transformation1} and \eqref{transformation2}, we have
 \begin{equation}\label{zxeta1}
\begin{split}
z(t)&=\bar{z}(t)+\textbf{z}(v(t),w)\\
x(t)&=\bar{x}(t)+\textbf{x}(v(t),w)\\
&=b(w)U_{d}^{-1}\xi(t)+\textbf{x}(v(t),w)\\
&=b(w)U_{d}^{-1}(\bar{\xi}(t)+\hat{\xi}(t))+\textbf{x}(v(t),w)\\
\eta(t)&=\bar{\eta}(t)+\theta(v(t),w)+C\bar{x}(t)\\
&=\bar{\eta}(t)+\theta(v(t),w)+Cb(w)U_{d}^{-1}(\bar{\xi}(t)+\hat{\xi}(t)).\\
\end{split}
\end{equation}
Note that $\textbf{z}(v(t),w)$, $\textbf{x}(v(t),w)$, $b(w)$ and $\theta(v(t),w)$ are all smooth functions, and the boundaries of the compact sets $\mathbb{V}$ and $\mathbb{W}$ are known. Then  $\textbf{z}(v(t),w)$, $\textbf{x}(v(t),w)$, $b(w)$ and $\theta(v(t),w)$ are all bounded for all $t\in[0,\infty)$. Then $z(t)$, $x(t)$ and $\eta(t)$ are all bounded for all $t\in[0,\infty)$.
Thus we conclude that $x_{c}(t)$ exists and is bounded for all $t\in[0,\infty)$, i.e., Property 1) of  Problem \ref{Problem1} is also satisfied.

The proof is thus completed.
\end{Proof}

\subsection{One Technical Lemma}\label{appendix2}
\begin{Lemma}\label{LemmaA1}
Consider the following system:
 \begin{equation}\label{zetasystem1}
\begin{split}
\dot{\zeta}_{1}&=\varphi_{1}(\zeta_{1},\zeta_{2},\psi,\mu)\\
\dot{\zeta}_{2}&=\varphi_{2}(\zeta_{1},\zeta_{2},\psi,\mu)+b(\mu)u\\
\end{split}
\end{equation}
where $\zeta_{1}\in\mathbb{R}^{n_{\zeta_{1}}}$, $\zeta_{2}\in\mathbb{R}$, $\psi\in\mathbb{R}^{n_{\psi}}$, $u\in\mathbb{R}$ $\mu\in\Omega\subset\mathbb{R}^{n_{\mu}}$ with $\Omega$ being some compact subset, $\varphi_{1}(\zeta_{1},\zeta_{2},\psi,\mu)$, $\varphi_{2}(\zeta_{1},\zeta_{2},\psi,\mu)$ and $b(\mu)$ are sufficiently smooth functions with $\varphi_{1}(0,0,0,\mu)=0$, $\varphi_{2}(0,0,0,\mu)=0$ and $b(\mu)>0$ for all $\mu\in\mathbb{R}^{n_{\mu}}$. Assume that there exists a $\mathcal{C}^{1}$ function $V_{1}(\zeta_{1})$, such that, for all $\zeta_{1},\zeta_{2}, \psi$, and all $\mu\in \Omega$,
$\underline{\gamma}_{1}(\|\zeta_{1}\|)\leq V_{1}(\zeta_{1})\leq \bar{\gamma}_{1}(\|\zeta_{1}\|)$ and
$\frac{\partial V_{1}(\zeta_{1})}{\partial \zeta_{1}}\varphi_{1}(\zeta_{1},\zeta_{2},\psi,\mu)\leq-\gamma_{1}(\|\zeta_{1}\|)+\phi_{1}(\zeta_{2})+\pi_{1}(\psi)$,
where $\underline{\gamma}_{1}(\cdot)$, $\bar{\gamma}_{1}(\cdot)$ and $\gamma_{1}(\cdot)$ are some known class $\mathcal{K}_{\infty}$ functions with $\gamma_{1}(\cdot)$ satisfying $\lim_{s\rightarrow0^{+}}\sup(s^{2}/\gamma_{1}(s))<\infty$,  $\phi_{1}(\cdot)$ and $\pi_{1}(\cdot)$ are some known smooth positive definite functions. Then there exists a control law
 \begin{equation}\label{u}
\begin{split}
u=-\rho(\zeta_{2})\zeta_{2}+\nu
\end{split}
\end{equation}
and a $\mathcal{C}^{1}$ function $V_{2}(\zeta)$, such that, for all $\mu\in\Omega$, and all $\zeta,\psi,\nu$
 \begin{equation}\label{V2zeta1}
\begin{split}
\underline{\gamma}_{2}(\|\zeta\|)\leq V_{2}(\zeta)\leq \bar{\gamma}_{2}(\|\zeta\|)
\end{split}
\end{equation}
 \begin{equation}\label{dotV2zeta1}
\begin{split}
\frac{\partial V_{2}(\zeta)}{\partial \zeta}\varphi(\zeta,\psi,\nu,\mu)\leq&-\gamma_{2}(\|\zeta\|)+\nu^{2}+\pi_{2}(\psi)\\
\end{split}
\end{equation}
where $\rho(\cdot)$ is some smooth positive function, $\varphi(\zeta,\psi,\nu,\mu)=\mbox{col}(\varphi_{1}(\zeta_{1},\zeta_{2},\psi,\mu),\varphi_{2}(\zeta_{1},\zeta_{2},\psi,\mu)+b(\mu)(-\rho(\zeta_{2})\zeta_{2}+\nu))$, $\nu\in\mathbb{R}$, $\zeta=\mbox{col}(\zeta_{1},\zeta_{2})$,  $\pi_{2}(\cdot)$ is some known smooth positive definite function, $\underline{\gamma}_{2}(\cdot)$, $\bar{\gamma}_{2}(\cdot)$ and $\gamma_{2}(\cdot)$ are some known class $\mathcal{K}_{\infty}$ functions with $\gamma_{2}(\cdot)$ satisfying $\lim_{s\rightarrow0^{+}}\sup(s^{2}/\gamma_{2}(s))<\infty$.
\end{Lemma}
\begin{Proof}
By applying the changing supply pair technique in \cite{Sontag1}, given any smooth function $\Delta_{1}(\zeta_{1})\geq0$, there exists a $\mathcal{C}^{1}$ function $\bar{V}_{1}(\zeta_{1})$, such that,
for all $\mu\in \Omega$, and all $\zeta_{1},\zeta_{2}, \psi$,
 \begin{equation}\label{bV1zeta1}
\begin{split}
\underline{\beta}_{1}(\|\zeta_{1}\|)\leq \bar{V}_{1}(\zeta_{1})\leq\bar{\beta}_{1}(\|\zeta_{1}\|)
\end{split}
\end{equation}
 \begin{equation}\label{dotbV1zeta1}
\begin{split}
\frac{\partial \bar{V}_{1}(\zeta_{1})}{\partial \zeta_{1}}\varphi_{1}(\zeta_{1},\zeta_{2},\psi,\mu)\!\leq&-\!\Delta_{1}(\zeta_{1})\|\zeta_{1}\|^{2}\!+\!\bar{\phi}_{1}(\zeta_{2})|\zeta_{2}|^{2}+\bar{\pi}_{1}(\psi)\|\psi\|^{2}
\end{split}
\end{equation}
where $\underline{\beta}_{1}(\cdot)$ and $\bar{\beta}_{1}(\cdot)$  are some known class $\mathcal{K}_{\infty}$ functions,  $\bar{\phi}_{1}(\cdot)$ and $\bar{\pi}_{1}(\cdot)$ are some known smooth positive functions.

Note that $\varphi_{2}(\zeta_{1},\zeta_{2},\psi,\mu)$ is sufficiently smooth with $\varphi_{2}(0,0,0,\mu)=0$ for any $\mu\in\mathbb{R}^{n_{\mu}}$.  Then, by Lemma 7.8 of \cite{Huang1}, there exist some known smooth functions $l_{2}(\zeta_{1})$, $\kappa_{2}(\zeta_{2})$ and $\omega_{2}(\psi)$ satisfying $l_{2}(0)=0$, $\kappa_{2}(0)=0$ and $\omega_{2}(0)=0$, such that, for all $\mu\in\Omega$, $\zeta_{1}\in\mathbb{R}^{n_{\zeta_{1}}}$, $\zeta_{2}\in\mathbb{R}$, $\psi\in\mathbb{R}^{n_{\psi}}$,
 \begin{equation}\label{varphi2}
\begin{split}
|\varphi_{2}(\zeta_{1},\zeta_{2},\psi,\mu)|\leq l_{2}(\zeta_{1})+\kappa_{2}(\zeta_{2})+\omega_{2}(\psi)
\end{split}
\end{equation}
which implies
 \begin{equation}\label{varphi3}
\begin{split}
|\varphi_{2}(\zeta_{1},\zeta_{2},\psi,\mu)|^{2}\leq&( l_{2}(\zeta_{1})+\kappa_{2}(\zeta_{2})+\omega_{2}(\psi))^{2}\\
\leq&3|l_{2}(\zeta_{1})|^{2}+3|\kappa_{2}(\zeta_{2})|^{2}+3|\omega_{2}(\psi)|^{2}.\\
\end{split}
\end{equation}
Since $l_{2}(0)=0$, $\kappa_{2}(0)=0$ and $\omega_{2}(0)=0$, there exist some smooth positive functions $\bar{l}_{2}(\zeta_{1})$, $\bar{\kappa}_{2}(\zeta_{2})$ and $\bar{\omega}_{2}(\psi)$, such that,
 \begin{equation}\label{varphi4}
\begin{split}
&3|l_{2}(\zeta_{1})|^{2}\leq\bar{l}_{2}(\zeta_{1})\|\zeta_{1}\|^{2}\\
&3|\kappa_{2}(\zeta_{2})|^{2}\leq\bar{\kappa}_{2}(\zeta_{2})|\zeta_{2}|^{2}\\
&3|\omega_{2}(\psi)|^{2}\leq\bar{\omega}_{2}(\psi)\|\psi\|^{2}.
\end{split}
\end{equation}
Combining \eqref{varphi2}, \eqref{varphi3} and \eqref{varphi4}, for all $\mu\in\Omega$, $\zeta_{1}\in\mathbb{R}^{n_{\zeta_{1}}}$, $\zeta_{2}\in\mathbb{R}$, $\psi\in\mathbb{R}^{n_{\psi}}$, we have
 \begin{equation}\label{varphi5}
\begin{split}
&|\varphi_{2}(\zeta_{1},\zeta_{2},\psi,\mu)|^{2}\leq \bar{l}_{2}(\zeta_{1})\|\zeta_{1}\|^{2}+\bar{\kappa}_{2}(\zeta_{2})|\zeta_{2}|^{2}+\bar{\omega}_{2}(\psi)\|\psi\|^{2}.
\end{split}
\end{equation}

In addition, since $b(\mu)>0$ for all $\mu\in\mathbb{R}^{n_{\mu}}$, there exist some real numbers $b_{m}$ and $b_{M}$ such that $0<b_{m}\leq b(\mu)\leq b_{M}<\infty$ for all $\mu\in\Omega$.

Let $W_{2}(\zeta_{2})=\frac{1}{2}\zeta_{2}^{2}$ and $\underline{\beta}_{2}(s)=\bar{\beta}_{2}(s)=\frac{1}{2}s^{2}$. Then $\underline{\beta}_{2}(|\zeta_{2}|)\leq W_{2}(\zeta_{2})\leq\bar{\beta}_{2}(|\zeta_{2}|)$ for any $\zeta_{2}$.

We further let $V_{2}(\zeta)=\bar{V}_{1}(\zeta_{1})+W_{2}(\zeta_{2})$. Then, by Lemma 11.3 of \cite{ChenZHuang2015}, we can choose some  class $\mathcal{K}_{\infty}$ functions $\underline{\gamma}_{2}(s)\leq\min\{\underline{\beta}_{1}(s/\sqrt{2}),\underline{\beta}_{2}(s/\sqrt{2})\}$ and $\bar{\gamma}_{2}(s)\geq\bar{\beta}_{1}(s)+\bar{\beta}_{2}(s)$ such that \eqref{V2zeta1} is satisfied.
According to \eqref{dotbV1zeta1} and \eqref{varphi5},  for all $\mu\in\Omega$ and all $\zeta,\psi,\nu$, we have
 \begin{equation}\label{dotV2zeta2}
\begin{split}
&\frac{\partial V_{2}(\zeta)}{\partial \zeta}\varphi(\zeta,\psi,\nu,\mu)\\
=&\frac{\partial \bar{V}_{1}(\zeta_{1})}{\partial \zeta_{1}}\varphi_{1}(\zeta_{1},\zeta_{2},\psi,\mu)+\zeta_{2}(\varphi_{2}(\zeta_{1},\zeta_{2},\psi,\mu)+b(\mu)(-\rho(\zeta_{2})\zeta_{2}+\nu))\\
\leq&-\Delta_{1}(\zeta_{1})\|\zeta_{1}\|^{2}\!+\!\bar{\phi}_{1}(\zeta_{2})|\zeta_{2}|^{2}+\bar{\pi}_{1}(\psi)\|\psi\|^{2}
+\frac{1}{4}|\zeta_{2}|^{2}+|\varphi_{2}(\zeta_{1},\zeta_{2},\psi,\mu)|^{2}\\
&-b_{m}\rho(\zeta_{2})|\zeta_{2}|^{2}+b_{M}|\zeta_{2}||\nu|\\
\leq&-\Delta_{1}(\zeta_{1})\|\zeta_{1}\|^{2}\!+\!\bar{\phi}_{1}(\zeta_{2})|\zeta_{2}|^{2}+\bar{\pi}_{1}(\psi)\|\psi\|^{2}+\frac{1}{4}|\zeta_{2}|^{2}+\bar{l}_{2}(\zeta_{1})\|\zeta_{1}\|^{2}\\
&+\bar{\kappa}_{2}(\zeta_{2})|\zeta_{2}|^{2}+\bar{\omega}_{2}(\psi)\|\psi\|^{2}-b_{m}\rho(\zeta_{2})|\zeta_{2}|^{2}+\frac{b_{M}^{2}}{4}|\zeta_{2}|^{2}|+\nu^{2}\\
=&-(\Delta_{1}(\zeta_{1})-\bar{l}_{2}(\zeta_{1}))\|\zeta_{1}\|^{2}-\bigg(b_{m}\rho(\zeta_{2})-\bar{\phi}_{1}(\zeta_{2})-\bar{\kappa}_{2}(\zeta_{2})-\frac{1+b_{M}^{2}}{4}\bigg)|\zeta_{2}|^{2}\\
&+\nu^{2}+(\bar{\pi}_{1}(\psi)+\bar{\omega}_{2}(\psi))\|\psi\|^{2}.\\
\end{split}
\end{equation}
Choose $\Delta_{1}(\zeta_{1})\geq \bar{l}_{2}(\zeta_{1})+1$, $\rho(\zeta_{2})\geq b_{m}^{-1}(\bar{\phi}_{1}(\zeta_{2})+\bar{\kappa}_{2}(\zeta_{2})+\frac{1+b_{M}^{2}}{4}+1)$, $\gamma_{2}(\|\zeta\|)=\|\zeta_{1}\|^{2}+|\zeta_{2}|^{2}$ and $\pi_{2}(\psi)\geq(\bar{\pi}_{1}(\psi)+\bar{\omega}_{2}(\psi))\|\psi\|^{2}$.  Then for all $\mu\in\Omega$, and all $\zeta,\psi,\nu$, we have
  \begin{equation}\label{dotV2zeta3}
\begin{split}
&\frac{\partial V_{2}(\zeta)}{\partial \zeta}\varphi(\zeta,\psi,\nu,\mu)\\
\leq&-\|\zeta_{1}\|^{2}-|\zeta_{2}|^{2}+\nu^{2}+(\bar{\pi}_{1}(\psi)+|\omega_{2}(\psi)|^{2})\|\psi\|^{2}\\
\leq&-\gamma_{2}(\|\zeta\|)+\nu^{2}+\pi_{2}(\psi).\\
\end{split}
\end{equation}

Thus the proof is completed.
\end{Proof}
\begin{Remark}
Lemma \ref{LemmaA1} can be viewed as an extension of Proposition 2.1 in \cite{LiuLu1}. 
\end{Remark}

\end{appendix}

\end{document}